\documentclass{elsart}

\usepackage{graphicx}
\usepackage{psfrag} 
\usepackage[final]{showkeys} 
\usepackage{amssymb}
\usepackage{amsmath} 
\usepackage{color}
\usepackage{algorithm}
\usepackage{algorithmic}
\journal{J. Comput. Phys.}
\newcommand{\Field}[1]{{\boldsymbol{\mathrm{#1}}}} 
\newcommand{\Pvec}[1]{{\vec{#1}}}

\newcommand{\curl}{\, \mathbf{curl}\, } 
 
\newcommand{\mdiv}{\, \mathbf{div}\, }

\newcommand{\rnum}{\mathbb{R}} 
\newcommand{\cnum}{\mathbb{C}}
\newcommand{\znum}{\mathbb{Z}} 

\renewcommand{\Re}{\mathrm{Re\,}}
\renewcommand{\Im}{\mathrm{Im\,}}
\newcommand{\dd}[1]{\mathrm{d}\,#1}
\newcommand{\ivec}[1]{\overleftarrow{#1}}

\begin{document}

\begin{frontmatter}
\title{Domain Decomposition Method for Maxwell's Equations: Scattering off Periodic Structures}
\author{Achim Sch\"adle\thanksref{matheon}},
\author{Lin Zschiedrich},
\author{Sven Burger},
\author{Roland Klose\thanksref{matheon}},
\author{Frank Schmidt}
\thanks[matheon]{Supported by the DFG Research
Center \textsc{Matheon} "Mathematics for key technologies" in
Berlin.}
\ead{schaedle@zib.de}
\ead{lin.zschiedrich@jcmwave.com}
\ead{sven.burger@jcmwave.com}
\ead{klose@zib.de}
\ead{frank.schmidt@zib.de}
\address{ZIB Berlin, Takustr.~7, D-14195 Berlin, Germany \\
JCMwave GmbH, Haarer Str.~14a, D-85640 Putzbrunn, Germany}

\begin{abstract}
  We present a domain decomposition approach for the computation of the
  electromagnetic field within periodic structures. We use a
  Schwarz method with transparent boundary conditions at the interfaces of the
  domains. Transparent boundary conditions are approximated by the perfectly 
  matched layer method (PML). To cope with Wood anomalies appearing in periodic 
  structures an adaptive strategy to determine optimal PML parameters is developed. \\ 
  We focus on the application to typical EUV lithography line masks.
  Light propagation within the multi-layer stack of the EUV mask is treated analytically. 
  This results in a drastic reduction of the computational costs and 
  allows for the simulation of next generation lithography masks 
  on a standard personal computer.
\end{abstract}
\end{frontmatter}

\section{Introduction}
\label{sec:intro}

\begin{figure}
  \begin{center}
  
    \psfrag{IlluminatingPlaneWave}[lc][lc][0.8][0]{$\vec{k}=(k_{1},k_{2},k_{3})$ 
      Illuminating Plane Wave}
    \psfrag{Air}[lc][lc][0.8][0]{Air}
    \psfrag{Line}[lc][lc][0.8][0]{Line}
    \psfrag{MultiLayerStack}[lc][lc][0.8][0]{Multi Layer Stack}
    \psfrag{Substrate}[lc][lc][0.8][0]{Substrate}
    \psfrag{Pitch}[lc][lc][0.8][0]{Period $a$}
    \psfrag{x}[lc][lc][0.8][0]{$x_{1}$}
    \psfrag{y}[lc][lc][0.8][0]{$x_{2}$}
    \psfrag{z}[lc][lc][0.8][0]{$x_{3}$}

    \includegraphics[width = 0.5\textwidth]{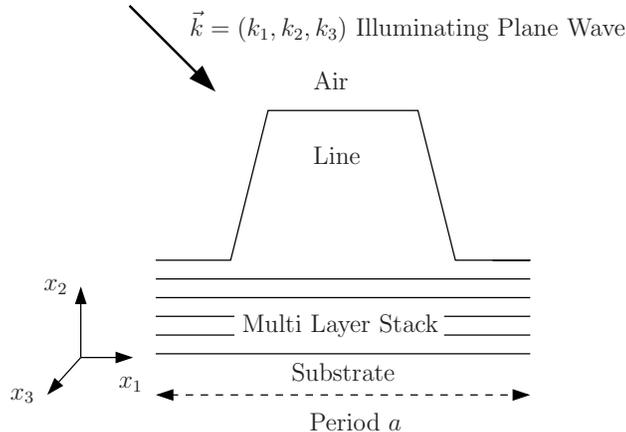}
    \caption{Layout of an EUV lithography line mask. The structure is periodically 
      repeated in $x_{1}$ direction and invariant in $x_{3}$ direction. The
      illuminating light is a plane wave with an arbitrary wave vector 
      $\vec{k} = (k_{1},k_{2}, k_{3}).$}
    \label{Fig.ScetchEUV}
  \end{center}
\end{figure}
The fabrication of semiconductor chips is essentially based on an optical projection
system. The pattern on a photolithography mask is transfered onto the chip by
optical projection. State of the art photolithography tools are operated with light
of a vacuum wavelength $\lambda \sim 193 \mathrm{nm},$~\cite{Burger2005bacus}.
Currently under development are tools that use extreme ultraviolet light (EUV) 
with a vacuum wavelength $\lambda \sim 13 \mathrm{nm}.$ 
These future systems will contain optical components including multi-layer structures 
which serve as mirrors. A typical section of an EUV lithography line mask is depicted in
Figure~\ref{Fig.ScetchEUV}. The line mask is invariant in $x_{3}$ direction and
periodic with period $a$ in $x_{1}$ direction. The multi-layer stack
may consist of more than $100$ layers. The thickness of each layer is below a
wavelength. The incident wave is twofold oblique -- oblique with respect to the mask
plane and oblique with respect to the multi-layer structure (conical incident).
The polarization of the incident field is arbitrary. 

In Section~\ref{sec:scattering} we introduce the mathematical setting of the 
arising scattering problem and derive the radiating boundary condition
in terms of Fourier modes. Further, we show that the exterior Dirichlet 
as well as Neumann boundary value problem is ill-posed in the presence of 
so called {\em Wood anomalies},~\cite{Petit1980a}. 
\\
To deal with the large computational domains we propose a domain decomposition
method (Section~\ref{sec:ddm}). 
The overall structure is split into sub-domains. The multi-layer sub-domain is
treated semi-analytically (c.f. Section~\ref{sec:tmm}) whereas the other subdomains 
are discretized by the finite element method utilizing the PML method 
to approximate transparent boundary conditions. The
PML method goes back to B\'erenger,~\cite{Berenger94}. Convergence
of the method was proven in~\cite{Lassas:98a,Lassas:01a} and \cite{Hohage03b} for
non-periodic problems. As is shown in Section~\ref{Sec:PML} the PML method
fails for periodic domains in the presence of Wood
anomalies. As a remedy we propose in Section~\ref{SubSec:AdaptedPML} 
a new automatic adaption of the layer size and the spatial discretization 
within the PML which leads to a quasi infinite layer thickness
in the presence of Wood anomalies. In Section~\ref{sec:variationalform} we 
introduce a variational formulation to couple
the PML to the interior problem.  
\\
In contrast to Elschner et al.~\cite{Elschner2002a,Elschner2000a} the
electromagnetic field $\Field{E} = (E_{1},E_{2},E_{3})$ 
is discretized with higher order Whitney elements 
for the $(E_{1},E_{2})$ component and Lagrange elements of the same order 
in the $E_{3}$ component. This allows for the
accurate evaluation of Fourier coefficients needed for the coupling to a multi-layer
stack as well for the computation of the far field coefficients.

The domain decomposition approach for the wave equation was first studied for the
scalar Helmholtz equation. Despr\'es and Shaidurov proposed to balance the energy 
fluxes across domain interfaces for Helmholtz problems, ~\cite{Despres90,Shaidurov91}. 
This idea was further expanded,~\cite{BenamouDepres97,Bourdonnaye97,CollinoGhanemiJoly00,Ghanemi98,CaiWidlund97,Ganderetal2002}. 
Toselli used the PML method at the interfaces of the sub-domains,~\cite{Toselli98}. 
This idea is closely related to the ideas of multiple
scattering, c.f.~\cite{Martin95} with further references. In each sub-domain a
simplified scattering problem is solved and the scattered field is added to the
incoming field for the neighboring domains. We have recently presented an additive
Schwarz algorithm for Helmholtz scattering problems with transparent boundary
conditions at the domain interfaces~\cite{SchZ05}. In this publication we used 
the $\mathrm{DtN}$ operator (Dirichlet to Neumann map). 
Since the definition of the $\mathrm{DtN}$ operator relies on the solvability of 
the exterior Dirichlet problem a $\mathrm{DtN}$ operator
may not exist for periodic structures. Hence we avoid the usage of the $\mathrm{DtN}$
operator for the formulation of the domain decomposition method in this paper.

The generalization to three dimensional geometries is straightforward.   

\section{Scattering off periodic line masks}
\label{sec:scattering}
The scattering off a periodic line mask is described by a Maxwell scattering problem,
with Bloch-periodic boundary condition in $x_{1}$ direction and transparent boundary
conditions in $x_{2}$ direction. The dependency on the $x_{3}$ component is
eliminated.

We consider electromagnetic scattering problems governed by the time-harmonic
Maxwell's equations
\begin{subequations}
  \label{eq.THMaxSystem}
  \begin{eqnarray}
    \curl \mu^{-1}\left(\Pvec{x} \right) \curl \Field{E}\left(\Pvec{x}\right) -
    \omega^{2} \varepsilon \left(\Pvec{x}\right)\Field{E}\left(\Pvec{x}\right) & = & 0,
    \label{eq.THMax}
    \\ \mdiv \varepsilon \left(\Pvec{x}\right) \Field{E}\left(\Pvec{x}\right) & = & 0,
    \label{eq.DivCond}
  \end{eqnarray}
\end{subequations}
which may be derived from the Maxwell's equations when assuming a
time dependence of the electric field as $\Field{E}(\Pvec{x}, t) =
\Field{E} \left(\Pvec{x} \right) \exp(-i \omega t)$ with angular frequency
$\omega$. The dielectric tensor $\varepsilon$ and the permeability tensor $\mu$ are
$L^{\infty}$ functions of the spatial variable $\Pvec{x}=(x_{1}, x_{2}, x_{3})$. In
addition we assume that the tensors $\varepsilon$ and $\mu$ do not depend
on $x_{3}$, that they are periodic functions in $x_{1}$ with period $a$,
i.e.~$\varepsilon(\Pvec{x}+(a, 0, 0)) = \varepsilon(\Pvec{x})$, $\mu(\Pvec{x}+(a, 0,
0)) =\mu(\Pvec{x})$, and that they are constant for $x_{2}>x_{2,+}$ and $x_{2}<x_{2,-}$ with
$x_{2,+}>x_{2,-}$. For simplicity assume that the dielectric and the permeability
tensors are isotropic so they may be treated as scalar valued functions. Recall that
any solution to~\eqref{eq.THMax} with $\omega \neq 0$ also meets the
divergence condition~\eqref{eq.DivCond}.
\\
A scattering problem may be defined as follows: given an incoming electric field
$\Field{E}_{\mathrm{inc}}$ satisfying the time-harmonic Maxwell's
equations~\eqref{eq.THMax} for $x_{2}>x_{2,+}$ and $x_{2}<x_{2,-}$, compute the total
electric field $\Field{E}$, which satisfies~\eqref{eq.THMax} in $\rnum^{3}$, such
that the scattered field $\Field{E}_{\mathrm{sc}} =
\Field{E}-\Field{E}_{\mathrm{inc}}$ defined for $x_{2}>x_{2,+}$ and $x_{2}<x_{2,-}$
meets the radiation condition given in Section~\ref{SubSec:HomoExtDom}. From a
physical point of view, the scattered field has to be outward radiating, so it only
transports energy towards infinity.
\\
It is possible to restrict the problem onto a two dimensional strip $[0, a] \times
\rnum$ provided that the incoming field is Bloch periodic in $x_{1}$,~\cite{Burger2004b} 
and depends
harmonically on $x_{3}$, i.e.
\begin{equation}
  \label{eq.BlochPeriodicity}
  \Field{E}_{\mathrm{inc}}\left(x_{1}+ a, x_{2}, x_{3} \right) = 
  \Field{\tilde{E}}_{\mathrm{inc}}
  \left(x_{1}, x_{2} \right) e^{i k_{1} a} e^{i k_{3} x_{3}}
\end{equation} 
where $\Field{\tilde{E}}_{\mathrm{inc}}$ is a periodic function in $x_{1}$ with
period $a$. The important case of an incoming plane wave meets these restrictions.
The total field $\Field{E}$ as well as the scattered field are then themselves Bloch
periodic in $x_{1}$ and depend harmonically on $x_{3}$.  This can be seen using a
symmetry argumentation where the unique solvability of the scattering problem is
assumed.

Below $\Field{E}$, $\Field{E}_{\mathrm{inc}}$ and
$\Field{E}_{\mathrm{sc}}$ denote the restriction of the respective field onto the
strip $[0, a] \times \rnum$. Let us introduce the domains $\Omega = [0, a] \times
[x_{2, -}, x_{2, +}]$, $\Omega_{+} = [0, a] \times [x_{2,+}, \infty]$ and accordingly
$\Omega_{-}$. With the definitions
\begin{eqnarray*}
  \curl_{3} \Field{E}  & = & 
  (\partial_{x_{2}}E_{3}-ik_{3} E_{2},\, ik_{3}E_{1}-\partial_{x_{1}}E_{3},\,
  \partial_{x_{1}}E_{2}-\partial_{x_{2}}E_{1})^{T}, 
  \\
  \mdiv_{3} \varepsilon \Field{E}  & = & 
  \partial_{x_{1}} \varepsilon E_{1}+
  \partial_{x_{2}} \varepsilon E_{2}+
  ik_{3} \varepsilon E_{3}
\end{eqnarray*} 
the scattering problem splits into an interior domain problem
\begin{gather}
  \label{eq.EInt}
  \begin{aligned}
    \curl_{3} \mu^{-1} \curl_{3} \Field{E}(x_{1}, x_{2}) - 
    \omega^{2} \varepsilon \Field{E}(x_{1}, x_{2})   = &  0 
    \quad (x_{1},x_{2}) \in \Omega,
    \\
    \Field{E}(0, x_{2})-\Field{E}(a, x_{2}) e^{ik_{1}a}  = & 0, 
  \end{aligned}
\end{gather}
an upper exterior domain problem 
\begin{gather}
  \label{eq.EExt}
  \begin{aligned}
    \curl_{3} \mu_{+}^{-1} \curl_{3} \Field{E}_{\mathrm{sc}, +}(x_{1}, x_{2}) - 
    \omega^{2} \varepsilon_{+} \Field{E}_{\mathrm{sc}, +}(x_{1}, x_{2})  = & 0 
    \quad (x_{1},x_{2}) \in \Omega_{+},
    \\
    \Field{E}_{\mathrm{sc},+}(0,x_{2})-\Field{E}_{\mathrm{sc},+}(a,x_{2}) e^{ik_{1}a}  
    = & 0 
  \end{aligned}
\end{gather}
and a lower exterior problem on $\Omega_{-}$ of similar type. 

Subproblems~\eqref{eq.EInt} and~\eqref{eq.EExt} are coupled by the following
matching conditions on the boundary $x_{2}=x_{2, +}$
\begin{subequations}
  \label{eq.Matching}
  \begin{eqnarray}
    \left(\Field{E} - 
      (\Field{E}_{\mathrm{sc}, +} +
      \Field{E}_{\mathrm{inc}, +}) 
    \right)
    \times \vec{n}_{+} & = & 0, \\
    \left(\mu^{-1} \curl_{3} \Field{E} - \left(
        \mu_{+}^{-1} \curl_{3} \Field{E}_{\mathrm{sc}, +} +
        \mu_{+}^{-1} \curl_{3}  \Field{E}_{\mathrm{inc}, +} \right) \right)
    \times \vec{n}_{+} &  = & 0,
  \end{eqnarray}
\end{subequations}
where $\vec{n}_{+} = (0,-1,0)^{T}$ denotes the unit normal vector.
An analogous condition holds on the boundary $x_{2}=x_{2,-}$, coupling
the interior and the lower exterior problem. 

\subsection{Radiation condition for homogeneous exterior domain problem}
\label{SubSec:HomoExtDom}
The exterior domain problems lack a radiation condition.  As both upper and lower
exterior problem can be treated similarly, we consider the upper exterior domain
problem only and drop the '$_{+}$'. Without loss of generality and to simply matters,
we assume that $x_{2, +}=0$.

Due to the periodicity of $\Field{E}_{\mathrm{sc}} \exp(-ik_{1}x_{1})$ 
the field has an expansion into Fourier modes
\begin{equation}
\label{Eqn:ExtFourierExpansion}
\Field{E}_{\mathrm{sc}}(x_{1}, x_{2}) = e^{+ik_{1}x_{1}} \sum_{n \in \znum}  
\vec{e}_{n}(x_{2}) e^{ix_{1} n 2\pi/a },
\end{equation}
with the Fourier coefficients
\[
\vec{e}_{n}(x_{2}) = 
\frac{1}{a} \int_{0}^{a} e^{-ik_{1}\xi} \Field{E}_{\mathrm{sc}}(\xi,x_{2}) 
e^{-i \xi (n2\pi/a)}\,d\xi~.
\]
The field $\Field{E}_n(x_{1}, x_{2}, x_{3}) = \vec{e}_{n}(x_{2}) \exp(i(n 2 \pi /
a+k_1)x_{1}) \exp(ik_{3}x_{3})$ is a solution of Maxwell's  equations~\eqref{eq.THMax}
for $x_{2}>0$. Hence inserting $\Field{E}_{n}$ in~\eqref{eq.THMax} yields
\begin{eqnarray*}
  \Field{E}_n \left(x_{1}, x_{2}, x_{3} \right) & = & 
  \vec{e}_{n, +} e^{i(n 2 \pi / a+k_{1})x_{1}} e^{ik_{2, n}x_{2}} e^{ik_{3}x_{3}} +
  \\ &&
   \vec{e}_{n, -} e^{i(n 2 \pi / a+k_{1})x_{1}} e^{-ik_{2, n}x_{2}} e^{ik_{3}x_{3}}, 
\end{eqnarray*}
with $k_{2, n} = \sqrt{k_0^2-(n 2 \pi / a + k_{1})^2 - k_{3}^2}$, where the branch
cut of the square root is along the negative real axis and $k_0 = \omega \sqrt{\mu
\varepsilon}$. From this representation it is easily seen that the field can be
decomposed into an incoming and an outgoing part.

We have to distinguish three cases:
\begin{enumerate}
\item $\Re{k_{2, n}} > 0$, $\Im{k_{2, n}}=0$ ({\em{propagating mode}})
  Both parts are propagating plane waves with wave vectors $(n 2 \pi / a + k_{1},
  k_{2, n}, k_{3})$ and $(n 2 \pi / a + k_{1}, -k_{2, n}, k_{3})$ respectively. The
  second part transports energy in the $-x_{2}$ direction.  We therefore require
  $\vec{e}_{n, -}=0$. This corresponds to the well known Sommerfeld radiation
  condition.
\item $\Re{k_{2, n}} = 0$, $\Im{k_{2, n}} > 0$ ({\em{evanescent mode}})
  The first part is evanescent in $x_{2}$ direction while the second part increases
  exponentially.  Therefore we again require $\vec{e}_{n, -}=0$.
\item $k_{2, n}=0$ ({\em{anomalous mode}})
  In this case both parts are equal and constant in $x_{2}$ direction. Energy is only
  transported in $x_{1}$ and $x_{3}$ directions. For the sake of a consistent
  notation we set $\vec{e}_{n, -}=0$.
\end{enumerate}  
Hence the correct radiation boundary condition is $\vec{e}_{n, -}=0$ for all $n \in
\znum$, such that the Fourier coefficients of the scattered field are given by
$\vec{e}_{n, \mathrm{sc}} = \vec{e}_{n, +}$. The anomalous case is rare. For example
for $k_{1}=0$ and $k_{3}=0$ it only occurs if $a=2\pi/(k_{0}n)$, hence the wavelength
must be a multiple of the period $a$.

\subsubsection{Ill-posed exterior Dirichlet/Neumann boundary value problems}
In our previous paper~\cite{SchZ05} the $\mathrm{DtN}$ operator was used to state the
coupling between the different domains. However the $\mathrm{DtN}$ operator must not
exist in the periodic setting -- the exterior Dirichlet problem is ill-posed in the
presence of anomalous modes.
 
This may be seen by rewriting Maxwell's equations separated in Fourier modes. With
$\vec{k}_{n} =(k_{1}+n2\pi/a, k_{2, n}, k_{3})$ the vectors $\vec{e}_{n,
\mathrm{sc}}$ satisfy the algebraic relations
\begin{subequations}
  \label{eq.THMaxModeSystem}
  \begin{eqnarray}
    - \vec{k}_{n} \times \left( \vec{k}_{n} \times \vec{e}_{n, \mathrm{sc}} \right) 
    - k_{0}^{2}\vec{e}_{n, \mathrm{sc}} & = & 0,
    \\
    \label{eq.THMaxMode}
    \vec{e}_{n, \mathrm{sc}} \cdot \vec{k}_{n} & = & 0.
    \label{eq.THMaxModeDiv}
  \end{eqnarray}
\end{subequations}
The first relation stems from Equation~\eqref{eq.THMax} and the second relation from
Equation~\eqref{eq.DivCond}. If $j \in \znum$ corresponds to an anomalous mode, that
is $\vec{k}_{j} = (k_{1}+j2\pi/a, 0, k_{3}),$ then $\vec{e}_{j, \mathrm{sc}} = (0, 1,
0)$ satisfies~\eqref{eq.THMaxModeSystem}. Hence $\Field{E}_{sc} = \vec{e}_{n,
\mathrm{sc}} \exp(i \vec{k_{j}} \cdot \vec{x})$ is a solution of the exterior domain
problem with zero Dirichlet tangential boundary values. Therefore the Dirichlet
boundary value problem is not uniquely solvable. Furthermore due to the divergence
condition~\eqref{eq.THMaxModeDiv} the vector $\vec{e}_{j, \mathrm{sc}}$ must be
perpendicular to $\vec{k}_{j}$. Hence for boundary values with $(d_{1}, 0,
d_{3})\cdot \vec{k}_{j} \neq 0$ the problem is not solvable at all.
\\
By an analogous argument one shows that the Neumann boundary value problem is also
ill-posed in the presence of anomalous modes.

\subsection{Scattering of an isotropic multi-layer stack -- the Transfer Matrix method}
\label{sec:tmm}
The Transfer Matrix method which according to~\cite{Plattner2003} was developed by
Schuster~\cite{Schuster49}, will be reviewed here shortly. 
For more details, the reader is referred
to~\cite{Plattner2003,Born99}.

Suppose we are in the situation of Figure~\ref{Fig.ScetchEUV}. Let us consider only
the material stack with $m$ finite layers, positioned at $x_{2,j}$, $j=0,\dots,m$,
with $x_{2,j}<x_{2,j+1}$. For $j=1,\dots,m$ the layer stack is given by the layer
thicknesses $x_{2,j}-x_{2,j-1}$, and the material coefficients $\epsilon_{j}$ and
$\mu_{j}$. Additionally for the semi-infinite half-spaces we have $\epsilon_{0}$,
$\mu_{0}$ and $\epsilon_{m+1}$, $\mu_{m+1}$. Since the Transfer Matrix algorithm is
applied to each Fourier mode $k_{1, n}$ separately we drop the subindex $n$ in this
section. In each layer we define the local wave vectors $\vec{k}_{j} = (k_{1},k_{2,
j},k_{3})$ and $\ivec{k}_{j} = (k_{1},-k_{2, j},k_{3})$ with $k_{2, j} =
\sqrt{\omega^{2}\epsilon_{j}\mu_{j}-k_{1}^{2}-k_{3}^{2}}$ such that $\Re{k_{2, j}}
\geq 0$ and $\Im{k_{2, j}} \geq 0.$ For a given excitation \footnote{Here it is not
assumed that the exciting field transports energy only in one direction.} from above
of the form $ \Field{E}_{\mathrm{inc}} = A_{m+1, \mathrm{inc}}
\exp(i\vec{k}_{m+1}\vec{x}) + B_{m+1,\mathrm{inc}} \exp(i\ivec{k}_{m+1}\vec{x}) $ we
want to calculate the reflected field $\Field{E}_{\mathrm{sc}} = A_{m+1,sc}
\exp(i\vec{k}_{m+1}\vec{x})$. From Snell's law we obtain that the field in each layer
is given by $\Field{E}_{j} = A_{j} \exp(i\vec{k}_{j}\vec{x}) + B_{j}
\exp(i\ivec{k}_{j}\vec{x})$. In the lower semi-infinite half space a purely outgoing
field is assumed, i.e $A_{0} =0$. In the layers we have $6m$ unknowns -- each $A$ or
$B$ has $3$ components. In the lower semi-infinite domain the only unknowns are the
three components of $B_{0}$ of the purely outgoing field. In the upper semi-infinite
domain there are six unknowns for the excitation and three for the reflected field.
\\
These unknown are determined by the following linear conditions arising from Maxwell's  
equations: there are $1 + 2 + 2m + 1$ equations from the divergence condition. 
At the $m+1$  boundaries of the layers there are $2(m+1)$
matching conditions for the tangential components of the Dirichlet data
and the same number of conditions from matching the Neumann data. 
$$
\left.
\begin{aligned}
\Field{E}_{j-1} \times \vec{n} &= \Field{E}_{j} \times \vec{n}
\\
\mu_{j-1} \curl \Field{E}_{j-1} \times \vec{n} &= \mu_{j} \curl\Field{E}_{j} \times \vec{n}  
& 
\end{aligned} \right\}
\mbox{ at } \vec{x} = \vec{x}_{j-1}
\mbox{ for } j = 1,\dots,n+1.
$$
Here $\Field{E}_{0} := \Field{E}_{\mathrm{inc}} + \Field{E}_{\mathrm{sc}}$.
The missing $4$ conditions are the tangential components of the Dirichlet and Neumann
data of the given incoming field.
\\
This yields a linear system of equations. To avoid large condition
numbers due to the complex material tensors, in each layer ansatz functions 
with amplitude equal to $1$ at the layer midpoint are used.

\section{Perfectly matched layer method}
\label{Sec:PML}
In the previous section we discussed the homogenous exterior domain problem and
derived transparent boundary conditions for each Fourier mode. Transforming back from
Fourier space, this boundary condition would be non-local and somehow the anomalous
case had to be treated separately. The perfectly matched layer method is an
approximate transparent boundary condition, introducing only small reflections that
are well under control. The reflections do not occur at the interface of the
computational domain and the PML, but stem from the truncation of the PML. Another
major advantage of the PML method is, that it fits in the finite-element framework,
described shortly in Section~\ref{SubSec:FEDiscr} and thus does not introduce
``full'' blocks in the discretization.
 
The PML method is based on a complex continuation of the scattered field. For
$\gamma=(1+i\sigma)$, $\sigma \geq 0$, we define the complex continued field
\begin{equation}
\Field{E}_{\gamma} = \sum_{n\in \znum}\vec{e}_{\mathrm{sc}, n} 
e^{i(n 2 \pi a+k_1)x_{1}} e^{ik_{n, 2} \gamma x_{2}} e^{ik_{3}x_{3}}.]
\label{eq:Egamma}
\end{equation}
With the definition 
\[
\curl_{3, \gamma} \Field{E}_{\gamma}  
=  
(\frac{1}{\gamma}\partial_{x_{2}}E_{\gamma, 3}-ik_{3} E_{\gamma, 2},
\, ik_{3}E_{\gamma, 1}-\partial_{x_{1}}E_{\gamma, 3},\,
\partial_{x_{1}}E_{\gamma, 2}-\frac{1}{\gamma}\partial_{x_{2}}E_{\gamma, 1})
\]
the field $\Field{E}_{\gamma}$ satisfies Maxwell's  equations~\eqref{eq.EExt} with
$\curl_{3}$ replaced by $\curl_{3, \gamma}$. In the absence of anomalous modes
$\Field{E}_{\gamma}$ is evanescent for $x_{2}\rightarrow \infty$,
\[
|\Field{E}_{\gamma} | \leq e^{-\kappa x_{2}} C,
\]
with $\kappa = \min_{n\in \znum}\{\Im{k_{n, 2}}, \sigma \Re{k_{n, 2}} \}$. The idea
is to restrict the complex continued exterior domain problem to a truncated domain
$\Omega_{\rho} = [0, a] \times [0, \rho]$ and to impose a zero Dirichlet boundary
condition at $x_{2}=\rho$. In case $\kappa$ is small or even $0$, i.e. if we are
``close'' to an anomalous mode a special adaptive PML is used, where the thickness
$\rho$ is increased like $1/\kappa$ and the discretization points are distributed
with an exponentially increasing mesh width guaranteeing an effective discretization,
c.f. Section~\ref{SubSec:AdaptedPML}. Thus the unbounded exterior
problem~\eqref{eq.EExt} is replaced by the truncated exterior domain problem
\begin{gather}
  \label{eq.EExtPML}
  \begin{aligned}
    \curl_{3} \mu_{+}^{-1} \curl_{3} \Field{E}_{\gamma, \rho}(x_{1}, x_{2}) - 
    \omega^{2} \varepsilon_{+} \Field{E}_{\gamma, \rho}(x_{1}, x_{2})  = & 0 
    \quad (x_{1},x_{2}) \in \Omega_{\rho},
    \\
    \Field{E}_{\gamma,\rho}(0, x_{2})-\Field{E}_{\gamma,\rho}(a,x_{2}) e^{ik_{1}a}  
    = & 0,  
    \\
    \Field{E}_{\gamma,\rho}|_{x_{2}=\rho} \times \vec{n} = & 0.
  \end{aligned}
\end{gather}
This modified truncated exterior problem is coupled to the interior problem using the
modified matching conditions, c.f.~\eqref{eq.Matching}
\begin{subequations}
  \label{eq.MatchingPML}
  \begin{eqnarray}
    \left(\Field{E} - 
      (\Field{E}_{\gamma, \rho} +
      \Field{E}_{\mathrm{inc}}) 
    \right)
    \times \vec{n} & = & 0, 
    \\
    \left(\varepsilon \curl_{3} \Field{E} - \left(
        \varepsilon_{+} \curl_{3, \gamma}  \Field{E}_{\gamma, \rho} +
        \varepsilon_{+} \curl_{3}  \Field{E}_{\mathrm{inc}} \right) \right)
    \times \vec{n} &  = & 0.
  \end{eqnarray}
\end{subequations}
Instead of homogenous Dirichlet boundary conditions at $x_{2}=\rho$ in
\eqref{eq.EExtPML} one can equally well require that the Neumann data is homogenous,
as in the PML the solution is oscillating and exponentially damped.

\subsection{Automatic Adaption of PML}
\label{SubSec:AdaptedPML}

\begin{algorithm}
  \caption{Adaptive PML method}
  \begin{algorithmic}[]
    \label{Algorithm:AdaptivePML}
    \REQUIRE $\epsilon, \sigma, h_{\mathrm{int}}, \kappa_{\mathrm{min}}$ 
    \STATE Compute $N_{\mathrm{p.w}}$ and $\xi_{\mathrm{max}}$ depending on $h_{\mathrm{int}}$ and 
    finite element order 
    \WHILE{(not converged)}
    \STATE 
    $\xi_{0}=0.0; \xi_{1}=h_{\mathrm{int}};N=1;$
    \WHILE {($-\ln(\epsilon)/(\xi_{N}\sigma)<\kappa_{\mathrm{min}}$)}
    \STATE
    $\xi_{N+1} = \xi_{N}+
    \max\{h_{\mathrm{int}},\; 2\pi \sigma \xi_{N} /(-\ln(\epsilon))/N_{\mathrm{p.w}} \}.
    $
    \IF {($\xi_{N+1} > 1/\epsilon$)}
    \STATE 
    \bf{break}
    \ELSE
    \STATE
    $N=N+1$
    \ENDIF
    \ENDWHILE
    \STATE Compute solution $u$ with PML discretization 
    $\{\xi_{0}, \xi_{1}, \dots, \xi_{N}\}$
    \IF {$\|u(\cdot, \xi_{N})\| \leq \epsilon \|u(\cdot)\|$}
    \STATE
    converged
    \ELSIF {$\xi_N>\xi_{\mathrm{max}}$}
    \STATE 
      break
    \ELSE
    \STATE
    $\kappa_{\mathrm{min}} = \kappa_{\mathrm{min}}/2$
    \ENDIF
    \ENDWHILE
  \end{algorithmic}
\end{algorithm}
As discussed in Section~\ref{SubSec:HomoExtDom} and Section~\ref{Sec:PML} the PML
method intrinsically fails in the presence of anomalous modes. For an anomalous mode
the field behaves like $\exp(i(k_{1}x_{1}+k_{3}x_{3}))$ and hence a complex
continuation in $x_{2}$ direction has no effect on the decay property of the field.
To obtain an effective transparent boundary condition we exploit the very specific
behavior in $x_{2}$ direction of the field and propose a mixed {\em a priori} and
{\em a posteriori} refinement strategy of the perfectly matched layer method
including the automatic adaption of the layer thickness $\rho.$ The algorithm we
propose is not restricted to the 2D periodic setting and was first published
in~\cite{Zschiedrich2006a} We therefore start from a
simple generic model. As in our paper~\cite{Schmidt02H,Zschiedrich03a} a prismatoidal
coordinate system in the exterior domain with a radial like coordinate $\xi$ and an
angular like variable $\eta$ is used. In the 2D periodic setting $\xi$ is simply the
$x_{2}$ coordinate and $\eta = x_{1}.$ Another example is a spherical coordinate
system in 3D $(r, \phi, \theta)$ with $\xi = r$ and $\eta = (\phi, \theta).$ To
proceed we assume the following expansion of the field in the exterior domain
\begin{equation}
  \label{Eqn:PMLExpansion}
  \Field{u}\left(\eta, \xi \right) \sim \int \Field{c}(\eta, \alpha) 
  e^{ik_{\xi}(\alpha)\xi}\,\dd{\alpha}
\end{equation} 
with $\Re k_{\xi}(\alpha) \geq 0$ and $\Im k_{\xi}(\alpha)\geq 0$. Hence in $\xi$
direction the field is a superposition of outgoing or evanescent plane waves. In the
periodic setting such an expansion is explicitly given
in~\eqref{Eqn:ExtFourierExpansion}.


The complex continuation, 
$\xi \mapsto \gamma \xi$ with $\gamma = 1+i\sigma$, gives
\begin{equation}
  \label{Eqn:PMLExpansionAbs}
  \|\Field{u}_{\gamma}\left(\eta, \xi \right)\| \sim 
  \int \|\Field{c}(\eta, \alpha)\| e^{-\kappa \xi}
\end{equation} 
with $\kappa = \sigma \Re k_{\xi} +\Im k_{\xi}.$
\\
The PML method only effects the outgoing part with $\Re k_{\xi}$ strictly larger
zero.  Field contributions with a large $\Re k_{\xi}$ component are efficiently
damped out.  Furthermore evanescent field contributions are damped out independently
of the complex continuation.  For a proper approximation of the oscillatory and
exponential behavior a discretization fine enough is needed to resolve the field.  In
contrast to that anomalous modes or ``near anomalous'' modes with $k_{\xi} \sim 0$
are neither evanescent nor damped out efficiently by the PML. Hence they 
enforce the usage of a large $\rho$ but can be well approximated with a relatively
coarse discretization in $\xi$ due to their smoothness in $\xi.$  
These requirements can only be satisfied by using an
adaptive discretization.  It is useful to think of the complex continuation as a
high-frequency filter. With a growing distance $\xi$ to the interior coupling
boundary the higher frequency contributions are damped out so that the discretization
can be coarsened.

For a given threshold $\epsilon$ selected according to the global accuracy requirements 
as described later we introduce the cut-off function
\[
\kappa_{\mathrm{co}, \epsilon}(\xi) = -\ln(\epsilon)/\xi~.
\]
With that at $\xi'>0$ each component in the expansion~(\ref{Eqn:PMLExpansionAbs}) with
$\kappa>\kappa_{\mathrm{co, \epsilon}}(\xi')$ is damped out by a factor smaller than the
threshold $\epsilon,$
\[
e^{-\kappa \xi'} < e^{-\kappa_{\mathrm{co, \epsilon}, \epsilon}(\xi') \xi} = 
e^{\ln(\epsilon)} = \epsilon.
\]
Assuming that this damping is sufficient we are allowed to select a discretization
which must only approximate the lower frequency parts with $\kappa \leq
\kappa_{\mathrm{co, \epsilon}}(\xi)$ for $\xi>\xi'.$ If we use a fixed number
$N_{\mathrm{p.w}}$ of discretization points per (generalized) wavelength $2\pi
/\kappa$ we get the following formula for the {\em a priori} determination of the
local mesh width $h({\xi}) = 2\pi \sigma
/\kappa_{\mathrm{co, \epsilon}}(\xi)/N_{\mathrm{p.w}}.$ A good choice of
$N_{\mathrm{p.w}}$ depends on the order of finite element used in $\xi-$direction and
need not to be adapted locally because the field depends smoothly in $\xi-$direction.
Since $\kappa_{\mathrm{co, \epsilon}}(\xi)
\rightarrow \infty$ for $\xi \rightarrow 0$ the local mesh width would be zero at $\xi=0.$
As it is not reasonable to use a finer discretization in the exterior domain than in
the interior domain we bound the local mesh width by the minimum mesh width
$h_{\mathrm{int}}$ of the interior domain discretization on the coupling boundary, 
\[
h({\xi}) = \max\{h_{\mathrm{int}},\; 2\pi \sigma /\kappa_{\mathrm{co, \epsilon}, 
  \epsilon}(\xi)/N_{\mathrm{p.w}} \}.
\] 
The parameters $\epsilon$ and $N_{\mathrm{p.w}}$ are also fixed accordingly to the
interior domain discretization quality. The grid $\{\xi_{0}, \xi_{1}, \xi_{2}, \dots
\}$ is recursively constructed by
\[
\xi_{n+1} = \xi_{n}+h(\xi_{n}).
\]
This way $\xi_{n}$ grows exponentially with $n.$ To truncate the grid we assume
that components in the expansion with $\kappa<\kappa_{\mathrm{min}}$
can be neglected so that the grid $\{\xi_{0}, \xi_{1}, \dots, \xi_{N}\}$ is
determined by $ \kappa_{\mathrm{co, \epsilon}, \epsilon}(\xi_{N}) <
\kappa_{\mathrm{min}}\leq\kappa_{\mathrm{co, \epsilon}, \epsilon}(\xi_{N-1}).$ In the
periodic setting there exists such a $\kappa_{\mathrm{min}}>0$ in case no anomalous
mode is present.
\\
As an {\em a posteriori} control we check if the field is indeed sufficiently
damped out at $\xi_{N}$, $\|u(\cdot, \xi_{N})\| \leq \epsilon \|u(\cdot)\|.$ \footnote{Here we
assume homogenous Neumann boundary conditions for the truncation of the
PML layer.  If homogenous Dirichlet boundary conditions are chosen for the truncation 
of the PML layer, the sufficient damping of the Neumann data may be checked instead.} 
Otherwise we recompute 
the solution with $\kappa_{\mathrm{min}} \rightarrow \kappa_{\mathrm{min}}/2$ 
\footnote{This strategy proved useful in many experiment. 
  However we consider to refine it.}. Since for an anomalous 
mode the field is not damped at all we restrict the maximum $\xi_N$ to
$\xi_N<\pi/k_0/\epsilon.$ The pseudocode to the algorithm is given
in Algorithm~\ref{Algorithm:AdaptivePML}.
\begin{figure}
  \begin{center}
    \psfrag{Einc}[lc][lc][1.2][0]{$\Field{E}_{\mathrm{inc}}$}
    \psfrag{Eout}[lc][lc][1.2][0]{$\Field{E}_{\mathrm{out}}$}
    \psfrag{Eref}[lc][lc][1.2][0]{$\Field{E}_{\mathrm{ref}}$}
    \includegraphics[width = 0.3\textwidth]{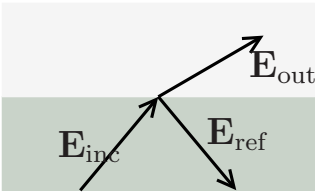}
    \caption{\label{Fig:TwoLayerScetch} Test problem for adaptive PML discretization. 
      The lower material has an refractive index equal to $n_{\mathrm{sub}}=1.5$, the
      upper material block consists of air ($n_{\mathrm{sup}}=1.0$). By Snell's law
      the field is totally reflected for an incident angle equal to the {\em critical
      angle} $\vartheta_{c}=180\cdot \mathrm{arcsin}(1.0/1.5)/\pi \approx 41.81.$
    }
  \end{center}
\end{figure}
\begin{figure}
  \begin{center}
    \includegraphics[width = 0.45\textwidth]{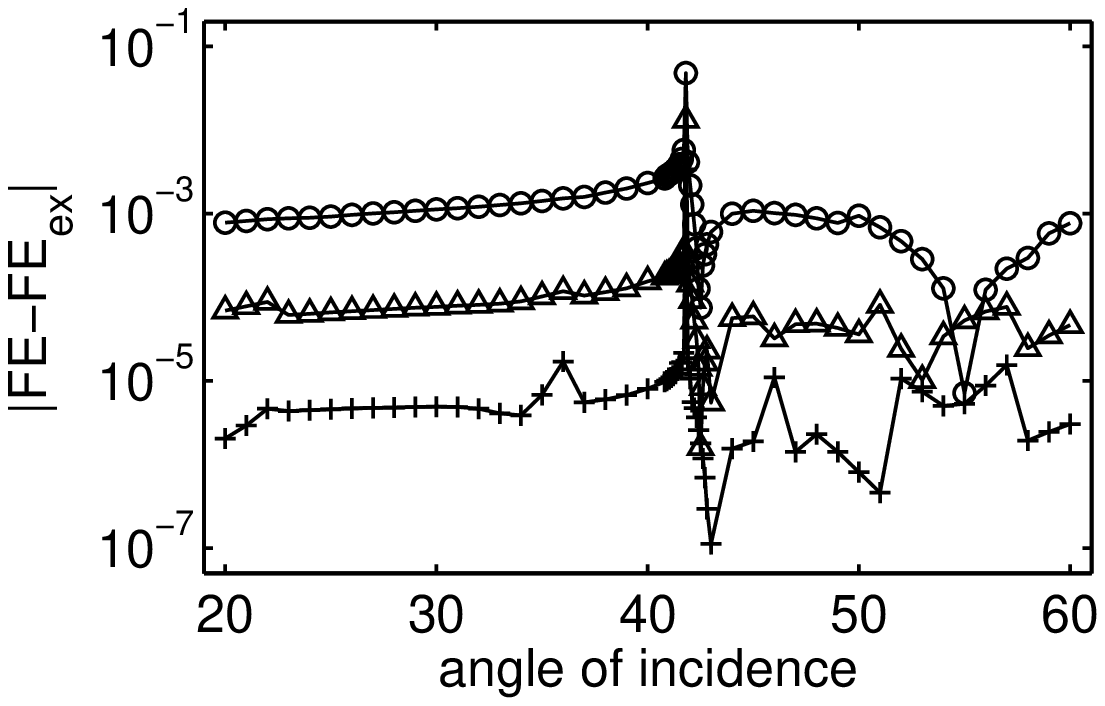}
    \includegraphics[width = 0.45\textwidth]{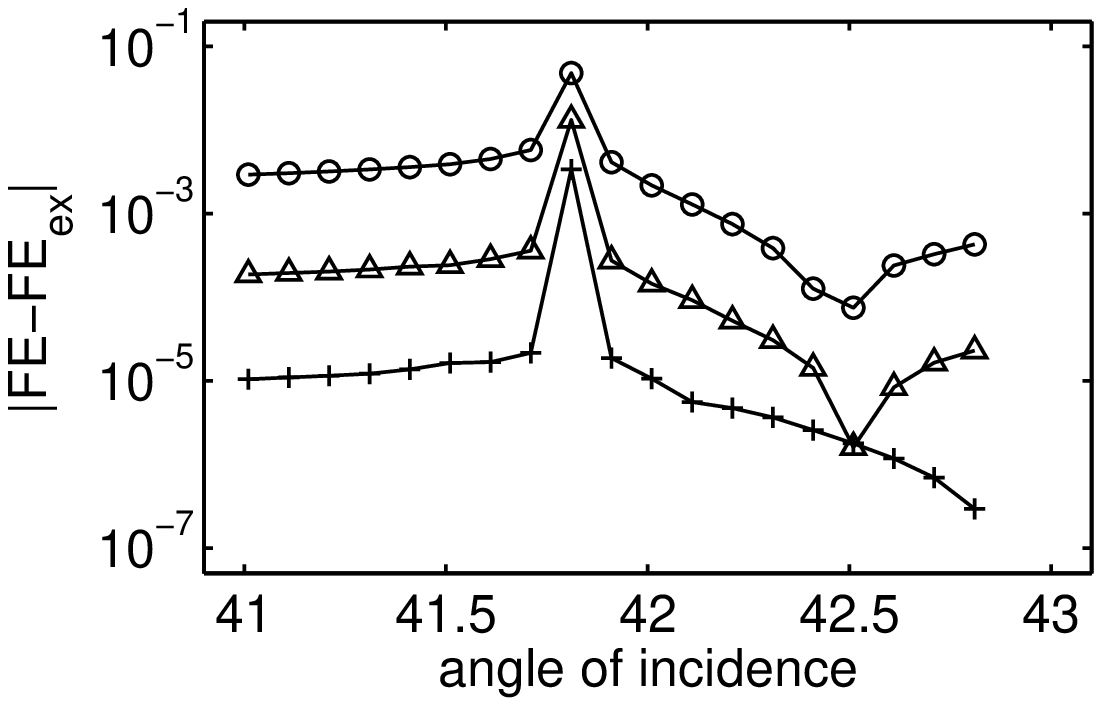}
    \caption{\label{Fig:PMLTestErrorScan} Left: 
      Field energy error in the interior domain. The three lines ($\circ$, $\hat{\quad}$,
      +) corresponds to different refinement levels of the interior domain.  Right: Zoom
      into left figure near critical angle.}
  \end{center}
\end{figure}

\begin{figure}
  \begin{center}
    \includegraphics[width = 0.45\textwidth]{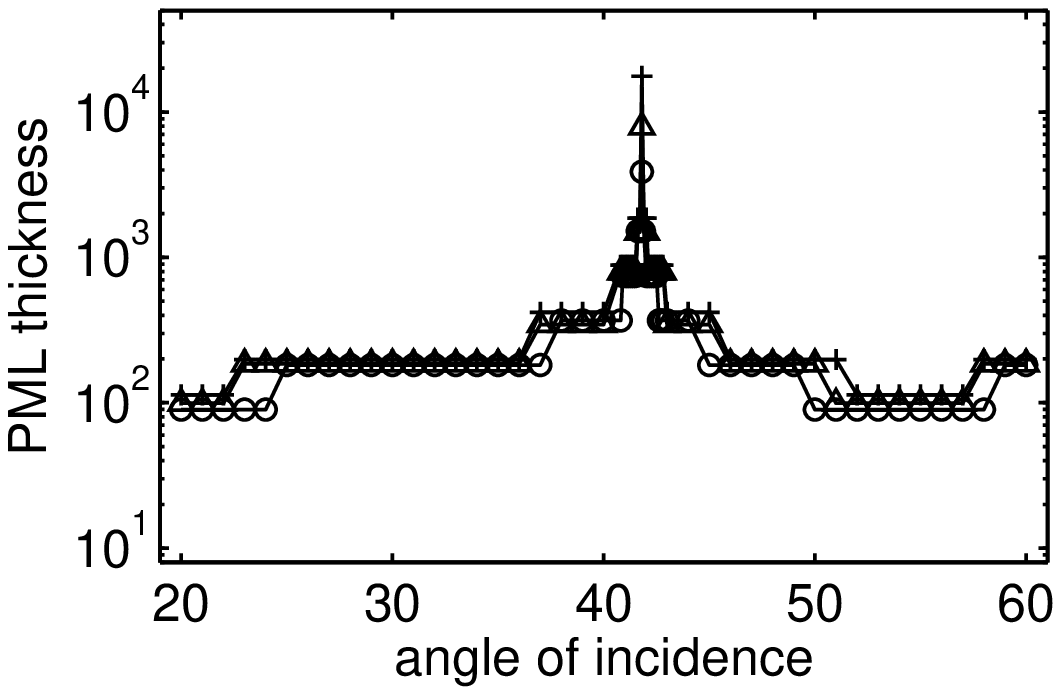}
    \includegraphics[width = 0.45\textwidth]{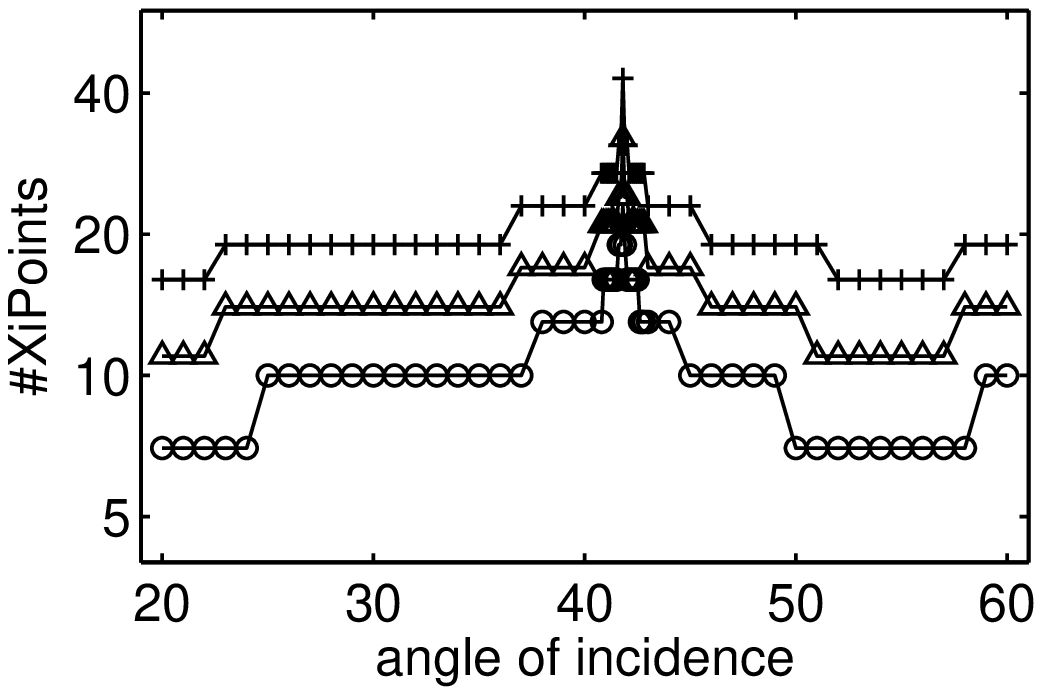}
    \caption{\label{Fig:PMLTestParameter} Left: 
      Thickness of the PML layer. At the critical angle the thickness is up to
      $10^{4}$ times larger than the diameter of the interior domain. Right: Number
      of discretization points $\xi_{j}$ used in the radial direction
      ($x_{2}$). Although the needed thickness of the layer is huge the number of
      unknowns used in the PML layer remains moderate.}
  \end{center}
\end{figure}
\begin{table}
  \begin{center}
    \begin{tabular}{rrr}
      Step & $\Delta E$ & $\Delta E'$  \\
      \hline
      0 & 0.359850   &  0.335129  \\
      1 & 0.159358   &  0.166207  \\ 
      2 & 0.048779   &  0.049502  \\
      3 & 0.012911   &  0.012912  \\
      4 & 0.003274   &  0.003266  \\
      5 & 0.000205   &  0.000820  \\
      6 & 0.000206   &  0.000205  \\
      7 & 0.000051   &  0.000051 
    \end{tabular}
  \end{center}
  \caption{Convergence of field energy at critical angle of incidence. The first
    column corresponds to the interior mesh refinement step. The relative error of
    the electric field energy in the interior domain is given in
    column two, $\Delta E = |\|\Field{E}_{\mathrm{ex}}\|_{L^{2}}^{2}-\|\Field{E}_{h}\|_{L^{2}}^{2}|/
    |\|\Field{E}_{\mathrm{ex}}\|_{L^{2}}^{2}.$ In column three the relative error of
    the magnetic field energy 
    $\Delta E' = |\|\curl \Field{E}_{\mathrm{ex}}\|_{L^{2}}^{2}-\|\curl \Field{E_{h}}\|_{L^{2}}^{2}|/
    |\|\curl \Field{E}_{\mathrm{ex}}\|$ is given. For fixed PML thickness the solution converges as 
    the interior mesh is refined.
    \label{Tab:ConvergencAtWA}}
\end{table}
To demonstrate the performance of the adaptive PML algorithm we compute the
reflection of a plane wave at a material jump,
c.f. Figure~\ref{Fig:TwoLayerScetch}. We vary the angle of incidence from
$\vartheta=20^{\circ}$ to $\vartheta=60^{\circ}.$ Further the incoming field is
rotated along the $x_{3}$ axis by an angle of $45^{\circ},$ so that the incidence is
twofold oblique (conical). Hence the unit direction of the incoming field is equal
to $\hat{k} = (\cos 45^{\circ} \sin \vartheta, \cos \vartheta, \sin 45^{\circ} \sin
\vartheta).$ The interior domain we used has as size of $1.5 \times 1$ in wavelength
scales. To measure the error we compute the field energy within the interior domain
and compare it with the analytic value. In Figure~\ref{Fig:PMLTestErrorScan} the
error is plotted for different refinement levels of the interior domain. The ``+''
line corresponds to the finest level. In Figure~\ref{Fig:PMLTestParameter} the
automatically adapted thickness of the PML is plotted (left) and the number of
discretization points $N$ in $\xi$ direction (right). As expected a huge layer is
used at the critical angle, whereas the total number of discretization points remains
moderate. As can be seen in Figure~\ref{Fig:PMLTestErrorScan} the maximum error
appears at the critical angle. From that one may suspect a failure of the automatic
PML adaption. But a closer analysis reveals that the chosen discretization in the PML
layer is sufficient as can be seen from Table~\ref{Tab:ConvergencAtWA}. Here the thickness
of the perfectly matched layer has been fixed and we further refined
the interior domain. This way we observe convergence to the true solution but the
convergence rate is halved at the critical angle. Hence the maximum error at the
critical angle comes from an insufficient interior discretization. We conjecture that
this is due to a dispersion effect. Since the wave is traveling along the $x_{1}$
direction it reenters the periodic domain leading to large ``path length''.

\section{Variational form}
\label{sec:variationalform}
The coupled problem given by~\eqref{eq.EExtPML},~\eqref{eq.EInt}
and~\eqref{eq.MatchingPML} can be casted into a variational problem on the Sobolev
space $H_{0,\rho}\left(\curl_{3}, \Omega \cup \Omega_{\rho}\right)$ of $H(\curl_{3})$
fields with generalized zero Dirichlet values at $x_{2}=\rho$. 
 
For a given 
test function $\Field{\Phi} \in H_{0,\rho} \left(\curl_{3},\Omega\cup\Omega_{\rho}\right)$ the
following identity holds true,
\begin{eqnarray}
  \lefteqn{\gamma \int_{\Omega_{\rho}} \overline{\Field{\Phi}} \cdot \curl_{3, \gamma} 
    \mu^{-1} \curl_{3, \gamma} \Field{E}_{\gamma} }\quad 
  && \nonumber
  \\
  &=&
  \gamma \int_{\Omega_{\rho}} \overline{\curl_{3, \gamma} \Field{\Phi}} \cdot 
  \mu^{-1} \curl_{3, \gamma} \Field{E}_{\gamma}
  -\int_{x_{2}=0} \overline{\Field{\Phi}} \cdot \mu^{-1} \curl_{3} 
    \Field{E}_{\mathrm{sc}}
  \times \vec{n}, \label{IntegrationByPartPML} 
\end{eqnarray}
where $\Field{E}_{\gamma}(x_{1},x_{2},x_{3}) = 
\Field{E}_{\mathrm{sc}}(x_{1},\gamma x_{2},x_{3})$, c.f.~\eqref{eq:Egamma}.
We first proof this identity for 
$\gamma \in \rnum \setminus \left\{ 0 \right\}.$ 
Using the non-euclidian coordinate change 
\[
\mathrm{T}^{-1}\,: \, \left(x_{1}, x_{2}, x_{3} \right)
\mapsto \left(x_{1}, \gamma^{-1} x_{2}, x_{3} \right)
\]
and applying the transformation
rules for differential forms, see~\cite{Zschiedrich2005a}, one gets
\begin{subequations}
  \label{TransformedIntegrals}
  \begin{eqnarray}
    \int_{\Omega_{\gamma \rho}} \overline{\Field{\Phi}^{*}} \cdot \curl_{3} 
    \mu^{-1} \curl_{3} \Field{E}_{\mathrm{sc}}  & = & 
    \int_{\Omega_{\rho}} \overline{\Field{\Phi}_{*}} \cdot \curl_{3} 
    \mu^{-1}_{*} \curl_{3} \Field{E}_{*} \quad \mbox{and} 
    \\
    \int_{\Omega_{\gamma \rho}} \overline{\curl_{3} \Field{\Phi}^{*}} \cdot 
    \mu^{-1} \curl_{3} \Field{E}_{\mathrm{sc}} & = & \int_{\Omega_{\rho}} \overline{\curl_{3} 
    \Field{\Phi}_{*}} \cdot 
    \mu^{-1}_{*} \curl_{3} \Field{E}_{*}
  \end{eqnarray}
\end{subequations}
with 
\begin{subequations}
  \label{TransformedQuantities}
  \begin{eqnarray}
    \mu_{*} &=&  |\mathrm{J}| \mathrm{J}^{-1} \mu \mathrm{J}^{-\mathrm{t}}
    \\
    \Field{E}_{*}(x_{1},x_{2},x_{3}) &=& \mathrm{J}^{\mathrm{t}}
    \Field{E}_{\mathrm{sc}}(x_{1}, \gamma x_{2}, x_{3})
    \\
    \Field{\Phi}_{*}(x_{1}, x_{2}, x_{3}) &=& 
    \mathrm{J}^{\mathrm{t}} \Field{\Phi}(x_{1}, x_{2}, x_{3}) = 
    \mathrm{J}^{\mathrm{t}} \Field{\Phi}^{*}(x_{1}, \gamma x_{2}, x_{3}).
  \end{eqnarray}
\end{subequations}
$\mathrm{J}=\mathrm{diag}(1, \gamma, 1)$ is the constant Jacobian of $ \mathrm{T}$.
Note that $\Field{\Phi}_{*}$, $\Field{E}_{*}$ are the pulled back fields to $\Field{\Phi}^{*}$ and
$\Field{E}_{\mathrm{sc}}$ in the sense of differential form calculus.
\\
We have 
\begin{subequations}
  \label{eq:insert}
  \begin{eqnarray}
    \gamma \int_{\Omega_{\rho}} \overline{\Field{\Phi}} \cdot \curl_{3, \gamma} 
    \mu^{-1} \curl_{3, \gamma} \Field{E}_{\gamma} 
    = 
    \int_{\Omega_{\rho}} \overline{\Field{\Phi}_{*}} \cdot \curl_{3} \mu^{-1}_{*}
    \curl_{3} \Field{E}_{*}
    \\
    \gamma \int_{\Omega_{\rho}} \overline{\curl_{3, \gamma} \Field{\Phi}} \cdot 
    \mu^{-1} \curl_{3, \gamma} \Field{E}_{\gamma}
    = 
    \int_{\Omega_{\rho}} \overline{\curl_{3} \Field{\Phi}_{*}} \cdot \mu^{-1}_{*} \curl_{3}
    \Field{E}_{*}
  \end{eqnarray}
\end{subequations}
which is verified by inserting~\eqref{TransformedQuantities} and using
$\curl_{3, \gamma} = |\mathrm{J}|^{-1} \mathrm{J} \curl_{3} \mathrm{J}^{\mathrm{t}}$.
\\
On the other hand integration by parts yields
\begin{gather}
  \begin{aligned}
    \lefteqn{\int_{\Omega_{\rho}} \overline{\Field{\Phi}_{*}} \cdot \curl_{3} \mu^{-1}_{*}
      \curl_{3} \Field{E}_{*} }\quad &&
    \\
    &=& 
    \int_{\Omega_{\rho}} \overline{\curl_{3} \Field{\Phi}_{*}} \cdot \mu^{-1}_{*} \curl_{3}
    \Field{E}_{*} 
    -\int_{x_{2}=0} \overline{\Field{\Phi}_{*}} \cdot \left( \mu^{-1}_{*} \curl_{3}
      \Field{E}_{*} \times \vec{n} \right)
  \end{aligned}
  \label{eq:partint}  
\end{gather}
and a respective equation for $\Field{E}_{\mathrm{sc}}$, $\mu$, and 
$\Field{\Phi}^{*}$ with the domain of integration $\Omega_{\gamma \rho}$. 
These together with equations in~(\ref{TransformedIntegrals}) give 
\begin{equation}
  \int_{x_{2}=0} \overline{\Field{\Phi}_{*}} \cdot \left( \mu^{-1}_{*} \curl_{3} 
    \Field{E}_{*} \times \vec{n} \right) = 
  \int_{x_{2}=0} \overline{\Field{\Phi}^{*}} \cdot \left( \mu^{-1} \curl_{3} 
    \Field{E}_{\mathrm{sc}}  \times \vec{n} \right)~.
  \label{eq:bd}  
\end{equation}
Using that $\Field{E}_{\gamma} = \mathrm{J}^{-\mathrm{t}} \Field{E}_{*}$ and using
that the tangential components of $\Field{\Phi}^{*}$ are equal to $\Field{\Phi}$, one
derives from~\eqref{eq:insert} and~\eqref{eq:bd} the desired
identity~\eqref{IntegrationByPartPML} for real $\gamma$. Since each term is a
holomorphic function in $\gamma$ the identity~\eqref{IntegrationByPartPML} holds true
for $\gamma \in \cnum \setminus \{ 0 \}$.

The coupled problem given by~\eqref{eq.EExtPML},~\eqref{eq.EInt}
and~\eqref{eq.MatchingPML} in weak form is given by
\begin{gather}
  \label{eq.varprePML}
  \begin{aligned}
    &\int_{\Omega} \overline{\curl_{3} \Field{\Phi}} \cdot \mu^{-1} \curl_{3} \Field{E}
    - \omega^{2} \overline{\Field{\Phi}} \cdot \varepsilon \Field{E}\;
    +
    \\
    &\gamma \int_{\Omega_{\rho}} \overline{\curl_{3, \gamma} \Field{\Phi}} \cdot 
    \mu^{-1} \curl_{3, \gamma} \Field{E}_{\gamma}
    - \omega^{2} \overline{\Field{\Phi}} \cdot \varepsilon \Field{E}_{\gamma}
    \\
    &= - \int_{x_{2}=0} \overline{\Field{\Phi}} \cdot \mu^{-1} (\curl_{3} \Field{E}  
    -  \curl_{3} \Field{E}_{\mathrm{sc}}) \times \vec{n} \;
  \end{aligned}
\end{gather}
Due to the Neumann coupling condition $\curl_{3} \Field{E}\times\vec{n} = \curl_{3}
\Field{E}_{\mathrm{sc}}\times\vec{n}+\curl_{3} \Field{E}_{\mathrm{inc}}\times\vec{n}$ the
boundary term is equal to $\int_{x_{2}=0} \Field{\Phi} \cdot \mu^{-1} \curl_{3}
\Field{E}_{\mathrm{inc}} \times \vec{n}$. This is not yet the basis for a Galerkin
ansatz in $H_{0,\rho}\left(\curl_{3}, \Omega \cup \Omega_{\rho}\right)$ as there is a
jump of the Dirichlet data across the boundary $x_{2} = 0$, precisely
$\Field{E}_{\gamma}+\Field{E}_{\mathrm{inc}} = \Field{E}|_{x_{2}=0}$. Let
$\Pi ({\Field{E}}_{\mathrm{inc}}\times \vec{n}) \in
H_{0,\rho}\left(\curl_{3},\Omega_{\rho}\right)$ denote an extension of a field with
tangential Dirichlet data equal to $\Field{E}_{\mathrm{inc}}\times \vec{n}$ at
$x_{2}=0$ to $H_{0,\rho}\left(\curl_{3},\Omega_{\rho}\right)$ and add this to
$\Field{E}_{\gamma}$ to obtain
\begin{gather}
  \label{eq.varPML}
  \begin{aligned}
    &\int_{\Omega} \overline{\curl_{3} \Field{\Phi}} \cdot \mu^{-1} \curl_{3} \Field{E}
    - \omega^{2} \overline{\Field{\Phi}} \cdot \varepsilon \Field{E}\;
    +
    \\
    &\gamma \int_{\Omega_{\rho}} \overline{\curl_{3, \gamma} \Field{\Phi}} \cdot 
    \mu^{-1} \curl_{3, \gamma} (\Field{E}_{\gamma}+ 
    \Pi({\Field{E}}_{\mathrm{inc}}\times \vec{n}))
    - \omega^{2} \overline{\Field{\Phi}} \cdot \varepsilon
    (\Field{E}_{\gamma}+\Pi({\Field{E}}_{\mathrm{inc}}\times \vec{n}))
    \\
    &= -\int_{x_{2}=0} \overline{\Field{\Phi}} \cdot \mu^{-1} (\curl_{3} \Field{E}_{\mathrm{inc}} \times \vec{n}) \;+
    \\  
    &\gamma \int_{\Omega_{\rho}} \overline{\curl_{3, \gamma} \Field{\Phi}} \cdot 
    \mu^{-1} \curl_{3, \gamma}\Pi({\Field{E}}_{\mathrm{inc}}\times \vec{n})
    - \omega^{2} \overline{\Field{\Phi}} \cdot \varepsilon 
    \Pi({\Field{E}}_{\mathrm{inc}}\times\vec{n})~.
  \end{aligned}
\end{gather}
This motivates the definition of the
composed field $\Field{u} \in H_{0,\rho}\left(\curl_{3}, \Omega \cup
\Omega_{\rho}\right)$ by $\Field{u}|_{\Omega} = \Field{E}$ and
$\Field{u}|_{\Omega_{\rho}} = \Field{E}_{\mathrm{sc}} +
\Pi(\Field{E_{inc}}\times\vec{n})$ and of the following bilinear form:
\begin{equation}
  \label{eq:a}
  a(\Field{\Phi},\Field{u}) := 
  a_{\Omega}(\Field{\Phi}|_{\Omega},\Field{u}|_{\Omega}) +
  a_{\Omega_{\rho}}(\Field{\Phi}|_{\Omega_{\rho}},\Field{u}|_{\Omega_{\rho}})
\end{equation}
with
\begin{equation}
  \label{eq:aO}
  a_{\Omega}(\Field{\Phi},\Field{u}) := 
  \int_{\Omega} \overline{\curl_{3} \Field{\Phi}} \cdot \mu^{-1} \curl_{3} \Field{u} 
  - \omega^{2} \overline{\Field{\Phi}} \cdot\varepsilon  \Field{u}~,
\end{equation}
\begin{equation}
  \label{eq:aOrho}
  a_{\Omega_{\rho}}(\Field{\Phi},\Field{u}) := 
  \gamma \int_{\Omega_{\rho}} \overline{\curl_{3, \gamma} \Field{\Phi}} \cdot 
    \mu^{-1} \curl_{3, \gamma} \Field{u} 
    - \omega^{2} \overline{\Field{\Phi}} \cdot \varepsilon \Field{u}~.
\end{equation}
With
\begin{equation}
  \label{eq:bG}
  b_{\Gamma}(\Field{\Phi},\Field{\Psi}) := 
  \int_{\Gamma} \overline{\Field{\Phi}} \cdot \mu^{-1} \Field{\Psi}~.
\end{equation}
we end up with the variational problem:
find  $\Field{u} \in H_{0,\rho}\left(\curl_{3}, \Omega \cup \Omega_{\rho}\right)$
such that for all 
$\Field{\Phi} \in H_{0,\rho}\left(\curl_{3}, \Omega \cup \Omega_{\rho}\right)$
\begin{equation}
  \label{eq:varPML2}
  a(\Field{\Phi},\Field{u}) =
  a_{\Omega_{\rho}}(\Field{\Phi},\Pi(\Field{E}_{\mathrm{inc}}\times\vec{n})) -
  b(\Field{\Phi}, \curl_{3} \Field{E}_{\mathrm{inc}}\times\vec{n})~. 
\end{equation}
Here we have avoided the definition of a DtN-operator. The total field is calculated
as the solution of a coupled system (computational domain coupled to the PML), where
the Dirichlet and Neumann data enter the equation on the ``right-hand side''. If
$\Field{u}$ is a solution of Maxwell's equations~\eqref{eq.THMaxSystem}, the
integration by parts identity can be rewritten using these bilinear forms as
\begin{equation}
  \label{eq:partinta}
  a_{\Omega}(\Field{\Phi},\Field{u}) - 
  b(\Field{\Phi}, -\curl_{3} \Field{u}\times\vec{n}) = 0
\end{equation}
This formula will be useful to represent the Neumann data. 
Note that in~\eqref{eq:partinta} $\vec{n}$ is the ``inward'' normal
with respect to $\Omega$.

\subsection{Finite element discretization}
\label{SubSec:FEDiscr}
To discretize the variational problem~\eqref{eq:varPML2} we use vectorial finite
elements on a triangular mesh in the interior domain and on a quadrilateral mesh in
the PML. The three sub-meshes -- lower PML mesh, interior domain mesh and upper PML
mesh -- fit non-overlapping. In the PML we use a rectangular mesh $[0, x_{1, 2},
\dots, a] \times [x_{2, +}, x_{2, +}+\xi_{1}, \dots, x_{2, +}+\xi_{N}]$ where
$\xi_{1}, \dots, \xi_{N}$ are determined as described in Section~\ref{SubSec:AdaptedPML}.
Since the Sobolev space $H_{0, \rho}({\curl_{3}, \Omega \cup \Omega_{\rho}})$ is
isomorphic to $H_{0, \rho}({\curl_{\mathrm{2D}}, \Omega \cup \Omega_{\rho}}) \times
H^{1}_{0, \rho}(\Omega \cup \Omega_{\rho})$ with the two dimensional $\curl$ operator
$\curl_{\mathrm{2D}} (u_{1}, u_{2}) = \partial_{x_{1}}u_{2}-\partial_{x_{2}}u_{1}$ we
use higher order Whitney elements to discretize the first and second component of
the electric field and standard Lagrange elements for the third field component of the
same order.  This finite element space is also used for waveguide mode computations,
c.f.~\cite{Schmidt00a} and the references therein.

\begin{figure}[htbp]
  \centering  \hfill
  \includegraphics[width=0.4\textwidth]{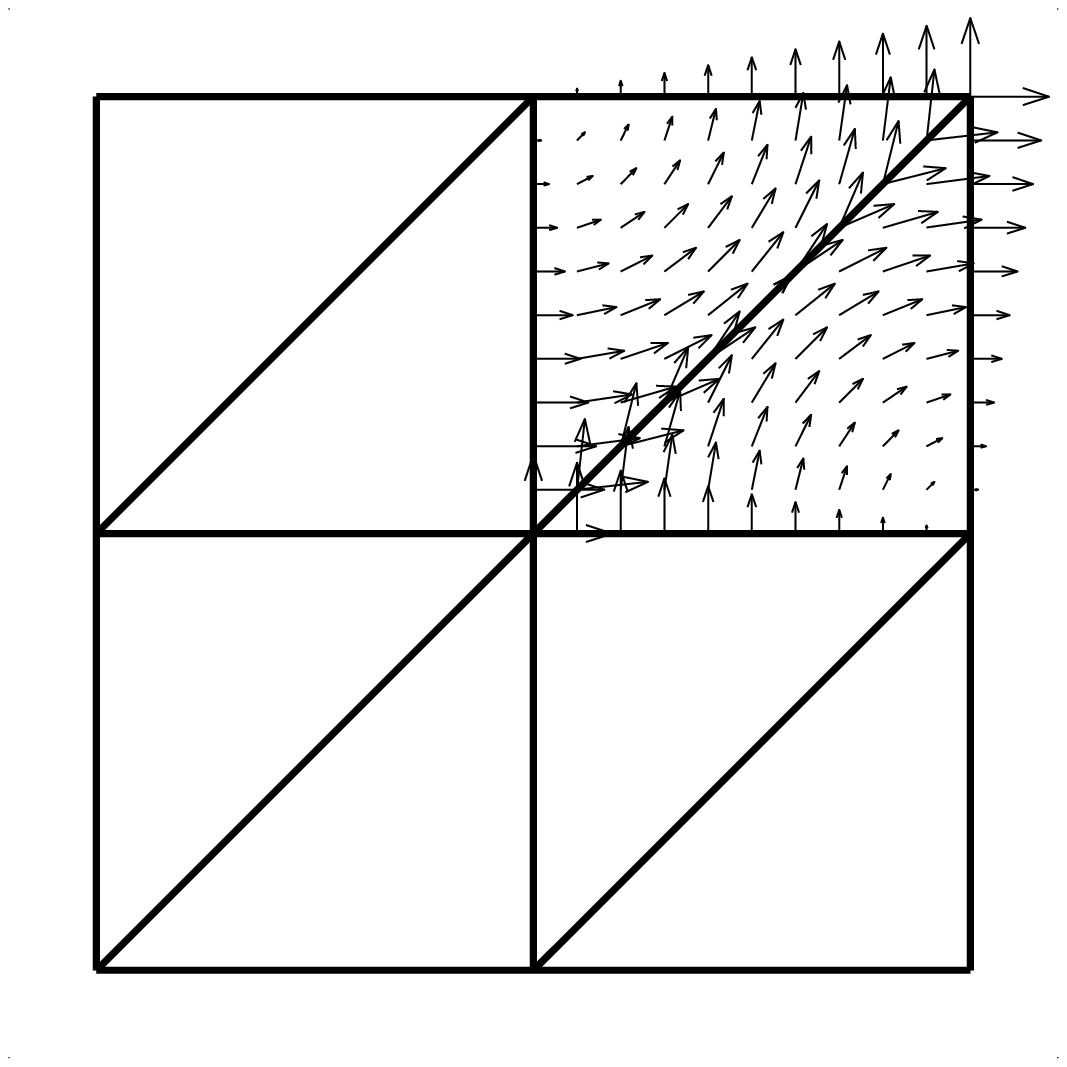}
  \hfill
  \includegraphics[width=0.4\textwidth]{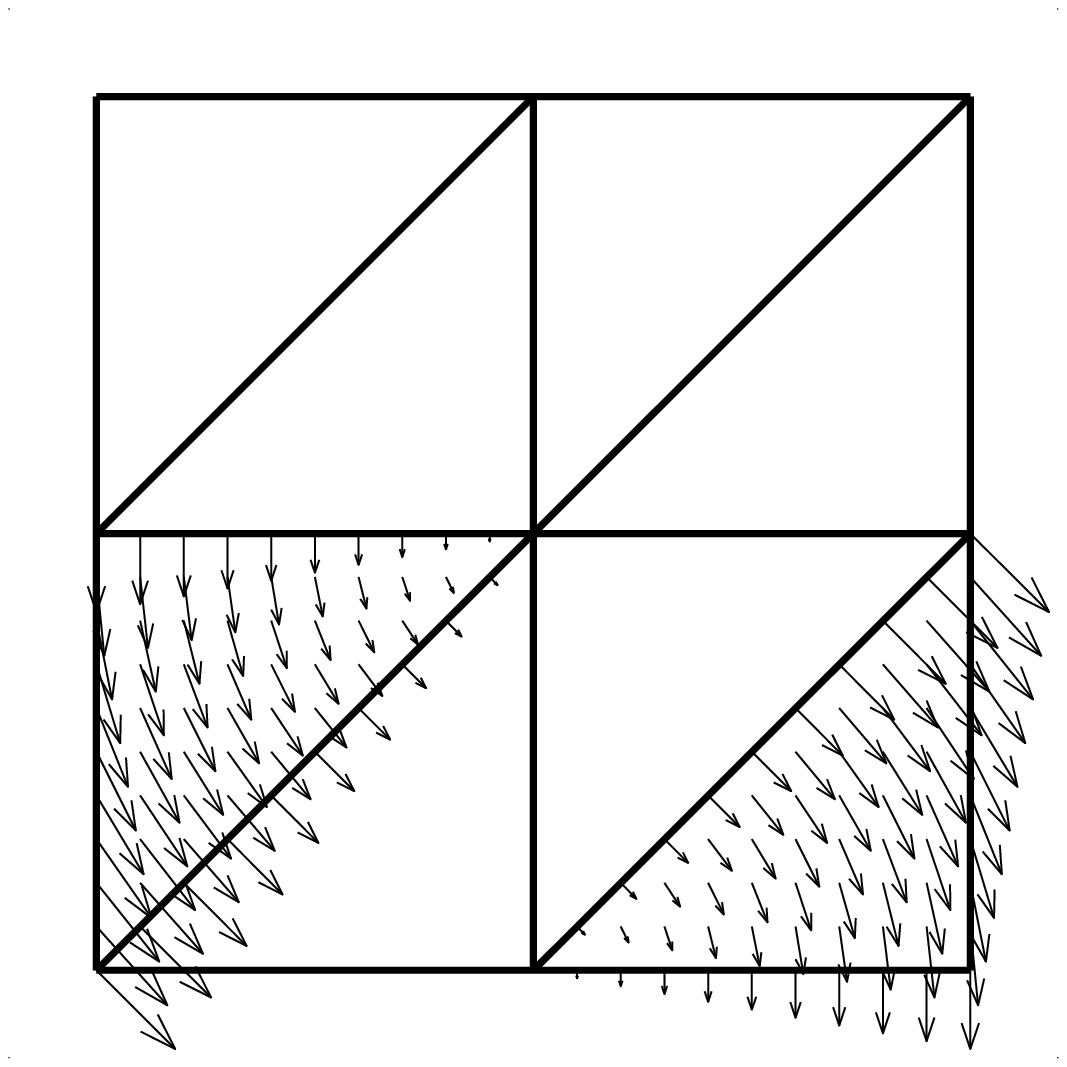}  \hfill
  \caption{First order edge elements on a simple grid. 
    In the interior the tangential component is continuous across 
    element boundaries. At the Bloch periodic boundary there is a phase
    shift.}
  \label{fig:edgeelements}
\end{figure}

Bloch periodicity is enforced by a multiplication of basis functions associated with
one of two corresponding periodic boundaries of the domain by the Bloch factor,
c.f.~\cite{Burger2004b}.  An interior edge element function remains unchanged, c.f.
Figure.~\ref{fig:edgeelements}~(left). The support of a basis function associated
with a periodic edge on the boundary consists of two triangles,
c.f. Figure.~\ref{fig:edgeelements}~(right). The restriction of the basis function to
the left triangle is defined as the standard shape function, whereas the
shape-function on the right triangle is multiplied by the Bloch factor $\exp(i
k_{1} a)$.  The construction of Bloch periodic Lagrange elements
is similar.

\section{Domain Decomposition Method}
\label{sec:ddm}
The idea for the Schwarz algorithm with transparent boundary conditions at the
interfaces is to calculate the solution on every sub-domain separately using
transparent boundary conditions and iteratively 
add the scattered field of each sub-domain to the
incoming field for the neighboring sub-domains. The presentation here
is restricted to the multiplicative Schwarz-algorithm.
\begin{figure}[tbp]
  \centering
  \psfrag{O1}{$\Omega_{1}$}
  \psfrag{O2}{$\Omega_{2}$}
  \psfrag{O3}{$\Omega_{3}$}
  \psfrag{O21}{$\Omega_{2,\rho,1}$}
  \psfrag{O23}{$\Omega_{2,\rho,3}$}
  \psfrag{O2ext}{$\Omega_{2,\rho}$}
  \includegraphics[width= 0.7\textwidth]{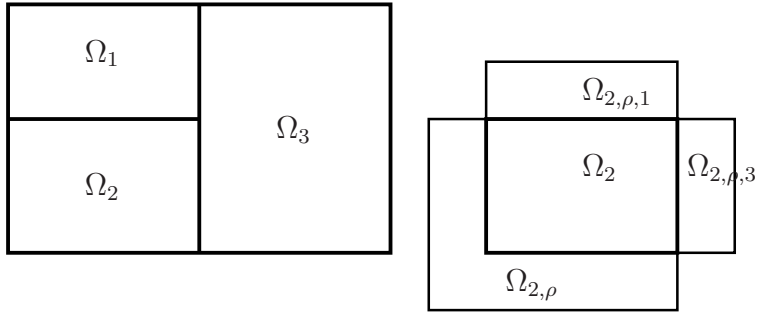}
  \caption{Schematic sketch of the various domains and PMLs. Left: The
  computational domain $\Omega$ is split in three sub-domains. Right:
  The sub-domain $\Omega_{2}$, with its three different PMLs.}
  \label{fig:scetch_dd}
\end{figure}
\\
In its most general form the domain-decomposition algorithm is given
in~\eqref{eq:ddalg}.
There $\Field{E}^{n}_{j}$ denotes the $n$th iterate
on sub-domain $\Omega_{j}$. $\Omega_{j,\rho,i}$ is
the PML domain to $\Omega_{j}$ at the interface to $\Omega_{i}$.
and by $\Omega_{j,\rho}$ we denote the PML domain to
$\Omega_{j}$ at the interface to the exterior, c.f. Figure~\ref{fig:scetch_dd}.
\begin{gather}
  \label{eq:ddalg}
  \begin{aligned}  
    \mbox{set}&\  \Field{E}_{j}=0 \mbox{ for all } j
    \\
    \mbox{wh}&\mbox{ile not converged}
    \\
    &\begin{aligned}
      \mbox{for}& \mbox{ all  sub-domains } j  
      \\
      &\begin{aligned}
        &\mbox{find $\Field{E}_{j}$ such that} 
        \\
        &\begin{aligned}
        a_{j}(\Field{\Phi},\Field{E}_{j}) 
        &=
        a_{\Omega_{j,\rho}}(\Field{\Phi},\Pi(\Field{E}_{\mathrm{inc}}\times\vec{n})) 
        + \sum_{i} 
        a_{\Omega_{j,\rho,i}}(\Field{\Phi},\Pi(\Field{E}_{i}\times\vec{n}))
        \\
        &
        -b_{\Gamma_{j}}(\Field{\Phi},\curl_{3}\Field{E}_{\mathrm{inc}}\times\vec{n}) 
        -\sum_{i} b_{\Gamma_{j,i}}(\Field{\Phi},\curl_{3}\Field{E}_{i}\times\vec{n})
        \end{aligned}
        \\
        &\forall \Field{\Phi} \in 
        H_{0,\rho}(\curl_{3},\Omega_{j} \cup \Omega_{j,\rho} \cup_{i} \Omega_{j,\rho,i} )
      \end{aligned} 
    \end{aligned}
  \end{aligned}
\end{gather}
This algorithm requires the evaluation of Neumann data
$\curl_{3}\Field{E}_{i}\times\vec{n}$ along the boundary. Moreover at cross
points, i.e at points given by $\Omega_{j,\rho,k}\cap\Omega_{j,\rho,i} \not =
\emptyset$ for $i\neq k$ or $\Omega_{j,\rho,k}\cap\Omega_{j,\rho} \not = \emptyset$, the
``incoming'' field is not a solution of Maxwell's equations, as it may have jumps.
\\
This difficulty maybe overcome by choosing sub-domains with a large overlap and by
coupling the incoming field to the computational not at the boundary but at some
additional artificial boundary in the interior. Here we do not pursue this strategy,
but avoid cross-points, by dividing the computational domain horizontally in several
sub-domains. Thus the sub-domains are arranged in a linear way. Each sub-domain has
only two well separated boundaries neglecting the periodic boundary, and at most two
neighboring domains.  Inserting an additional post-processing step, the Neumann-data
can be evaluated weakly.
\begin{gather}
  \label{eq:ddalweak}
  \begin{aligned}  
    \mbox{set}&\  \Field{E}_{j}=0 \mbox{ for all } j
    \\
    \mbox{wh}&\mbox{ile not converged}
    \\
    &\begin{aligned}
      \mbox{for}& \mbox{ all  sub-domains } j  
      \\
      &\begin{aligned}
        &\mbox{find $\Field{E}_{j}$ such that} 
        \\
        &\begin{aligned}
        a_{j}(\Field{\Phi},\Field{E}_{j}) 
        &=
        +a_{\Omega_{j,\rho}}(\Field{\Phi},\Pi(\Field{E}_{\mathrm{inc}}\times\vec{n})) 
        + \sum_{i} 
        a_{\Omega_{j,\rho,i}}(\Field{\Phi},\Pi(\Field{E}_{i}\times\vec{n}))
        \\
        &
        -b_{\Gamma_{j}}(\Field{\Phi},\curl_{3}\Field{E}_{\mathrm{inc}}\times\vec{n}) 
        -\sum_{i} b_{\Gamma_{j,i}}(\Field{\Phi},\curl_{3}\Field{E}_{i}\times\vec{n})
        \end{aligned}
        \\
        &\forall \Field{\Phi} \in 
        H_{0,\rho}(\curl_{3},\Omega_{j} \cup \Omega_{j,\rho} \cup_{i} \Omega_{j,\rho,i})
      \end{aligned} 
      \\
      \mbox{for}&\mbox{ all subdomains} j 
      \\
      &
      \sum_{i} b_{\Gamma_{j,i}}(\Field{\Phi},\curl_{3}\Field{E}_{j}\times\vec{n})
      +  b_{\Gamma_{j}}(\Field{\Phi},\curl_{3}\Field{E}_{j}\times\vec{n})
      =- a_{\Omega_{j}}(\Field{\Phi},\Field{E}_{j})
    \end{aligned}
  \end{aligned}
\end{gather}
In order to distinguish for $k
\not = l$ the contribution
$b_{\Gamma_{j,k}}(\Field{\Phi},\curl_{3}\Field{E}_{j}\times\vec{n})$ from, say
$b_{\Gamma_{j,l}}(\Field{\Phi},\curl_{3}\Field{E}_{j}\times\vec{n})$, 
it is required that there are no test functions that have a support in
elements adjacent to $\Gamma_{j,k}$ and $\Gamma_{j,l}$ simultaneously.

\subsection{Schwarz algorithm for EUV}
\label{SubSec:SAEUV}
For the special application -- scattering off an EUV-line mask -- one can
make use of the ``simple'' geometry of the double layer stack that serves as
a mirror employing the Transfer Matrix algorithm of Section~\ref{sec:tmm}. 

\begin{figure}
  \begin{center}
    \psfrag{y0}[lc][lc][0.8][0]{$x_{2,0}$}
    \psfrag{yh}[lc][lc][0.8][0]{}
    \psfrag{CD}[lc][lc][0.8][0]{Computational domain}
    \psfrag{M1}[lc][lc][0.8][0]{Material 0}
    \psfrag{OW}[lc][lc][0.8][0]{Outgoing wave}
    \psfrag{RW}[lc][lc][0.8][0]{Reflected wave}
    \psfrag{MLS}[lc][lc][0.8][0]{Multi-layer stack}
    \psfrag{S}[lc][lc][0.8][0]{Substrate}
        
    \includegraphics[width = 0.43\textwidth]{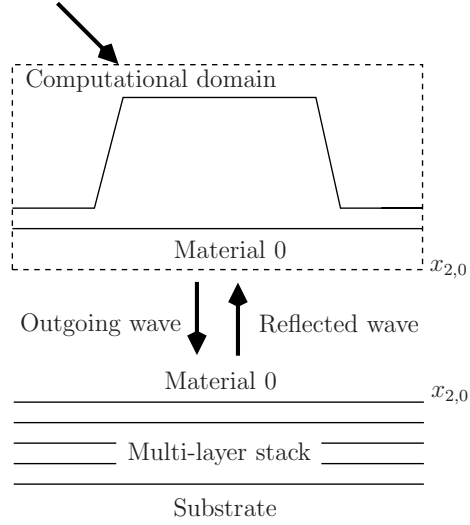}
    \caption{Decomposition of the problem into two infinite sub-domains. 
      The scattering problem is solved by the Finite Element Method in the upper
      domain and quasi-analytically in the lower domain.}
    \label{Fig:DomainDecompositionScetch}
  \end{center}
\end{figure}

A simple situation is depicted in Figure~\ref{Fig:DomainDecompositionScetch}.
The upper domain contains the mask line whereas the lower domain consists of the
multi-layer stack and the lower substrate block. Instead of solving Maxwell's
equations by the finite element method in the multi-layer stack, the incoming field
is Fourier transformed and for each Fourier mode the Transfer Matrix algorithm is
used to calculate the scattered field. This can even be simplified. If the tangential
component of each Fourier-mode vector field is written as the linear combination of
two linear independent polarizations, it is sufficient to compute the reflection
coefficients of the multi-layer layer stack for each mode and each polarization only
once. The number of Fourier mode ranges from $=n_{min}$ to $n_{max}$.  To determine
these, we set $k_{max} = 0.1\cdot 2\pi/h_{max}$, where $h_{max}$ is the maximum
segment size of a finite element at the boundary. Then $n_{max}$ is the greatest
integer such that $k_{1} + n_{max} 2\pi/a < k_{max}$, i.e and $n_{min}$ is the
greatest integer, such that $k_{x}-n_{min} 2\pi/a > -k_{max}$.

\section{Numerical Examples}
\label{Sec:NumEx}

\subsection*{An academic example}
\begin{figure}[tbp]
  \centering
  \psfrag{n1}{$\varepsilon_{8}$}
  \psfrag{n2}{$\varepsilon_{7}$}
  \psfrag{n3}{$\varepsilon_{6}$}
  \psfrag{n4}{$\varepsilon_{5}$}
  \psfrag{n5}{$\varepsilon_{3}$}
  \psfrag{n6}{$\varepsilon_{1}$}
  \psfrag{n7}{$\varepsilon_{4}$}
  \psfrag{n8}{$\varepsilon_{2}$}
  \psfrag{D1}{$\Omega_{1}$}
  \psfrag{D2}{$\Omega_{2}$}
  \includegraphics[height=2.5in]{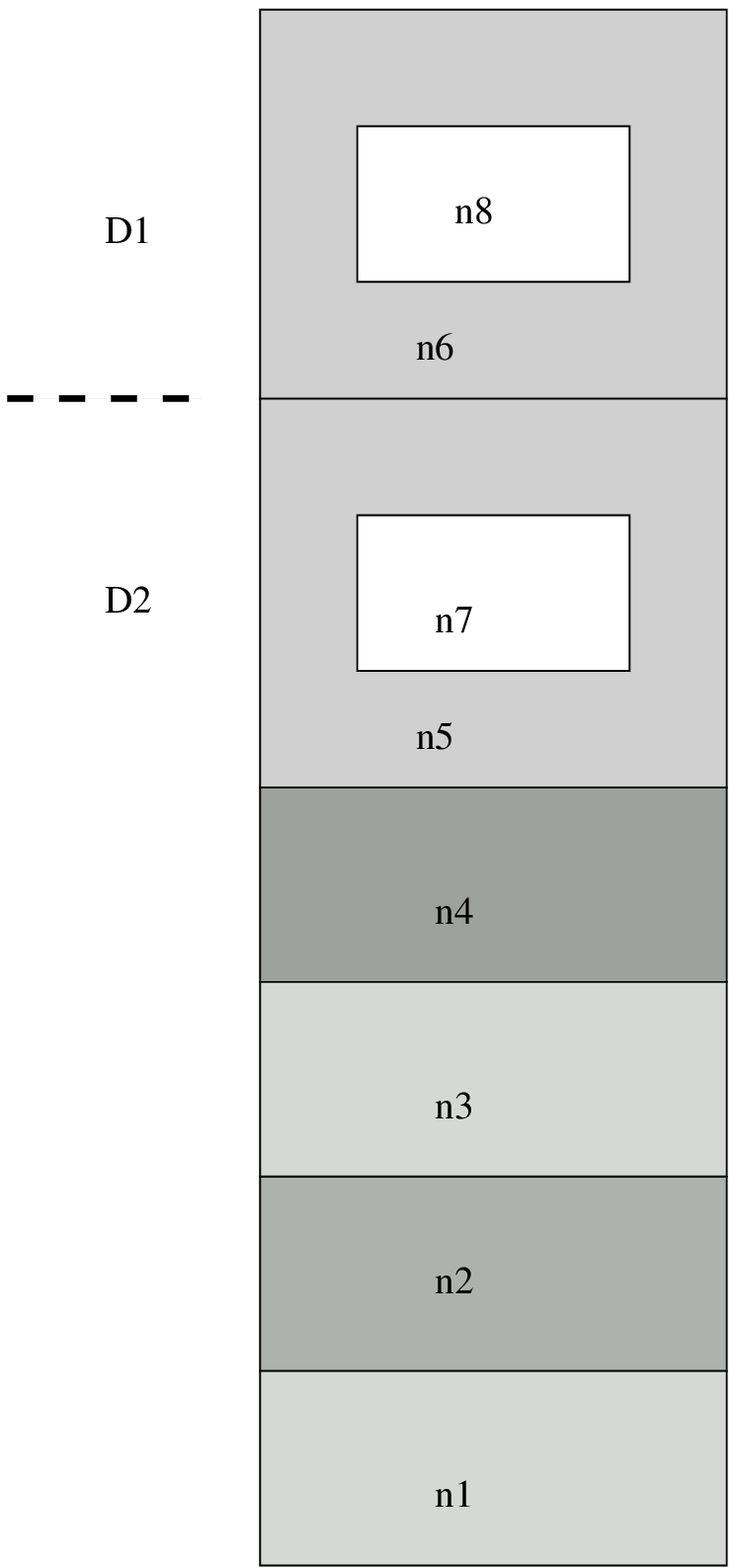}
  \hspace*{1ex}
  \includegraphics[height=2.5in]{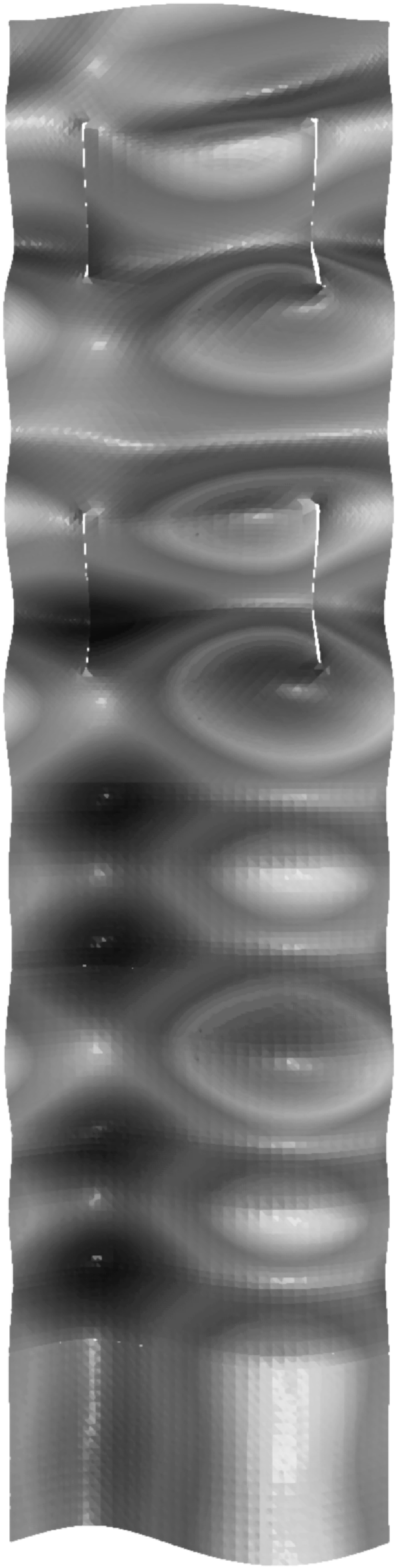}
  \hspace*{1ex}
  \includegraphics[height=2.5in]{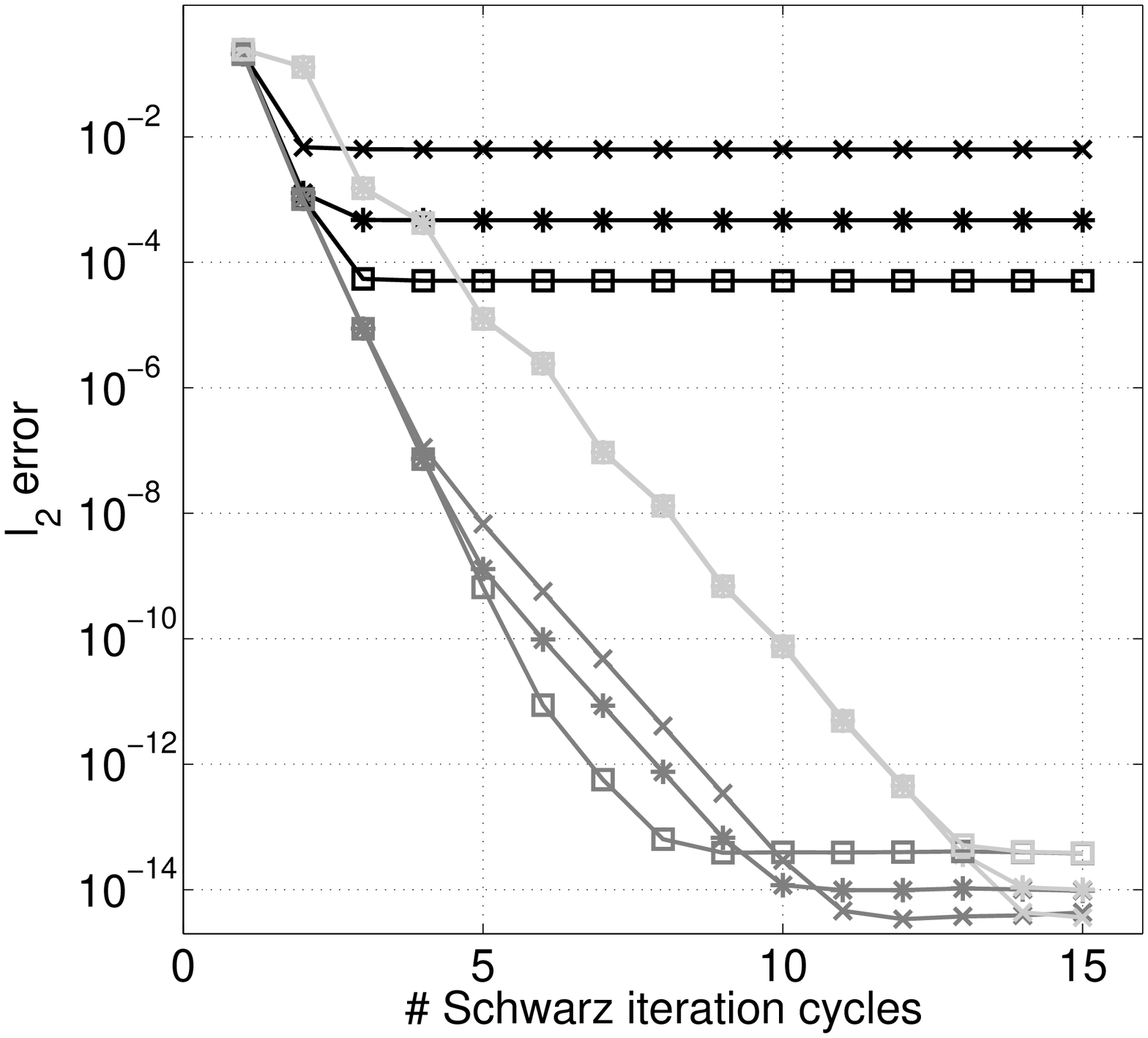}
  \caption{Material distribution (left) and magnitude of the electric field (middle)
    for a simple test problem. Convergence plot (right) for 
    $\vec{k}_{\mathrm{inc}}= (1, -2, 1)$, $
    \lambda = 0.84$ and different refinement levels.}
  \label{fig:example1}
\end{figure}

The simple geometry of this example is depicted in Figure~\ref{fig:example1}. It
consists of three domains $\Omega_{1}$ and $\Omega_{2}$, each with a quadrilateral
material inhomogeneity and $\Omega_{3}$ a layer stack of four layers below
$\Omega_{2}$. The period is $a=1$. 
The different shadings correspond to different materials as
indicated. The permeability is equal to $1$ everywhere. The permittivity is given by
$\varepsilon_{1}=1.01$, $\varepsilon_{2}=1.52$, $\varepsilon_{3}=1.03$,
$\varepsilon_{4}=1.54$, $\varepsilon_{5}=1.55$, $\varepsilon_{6}=1.06$,
$\varepsilon_{7}=1.57$, $\varepsilon_{8}=1.08$. The semi-infinite top and lower
strips, with refraction indices $\varepsilon_{9}=1$ and $\varepsilon_{0}=1$ are not
shown. These are completely modeled by the PML method. Hence there are rather big
jumps in the material coefficients at domain interfaces.

The incoming field is a plane waves with wave vector $\vec{k}_{\mathrm{inc}}= (1,-2,1)$ 
and wave length, $\lambda = 0.84$.  The strength is $\vec{s}_{\mathrm{inc}}=
(1,1,1)\times\vec{k_{\mathrm{inc}}}/||(1,1,1)\times\vec{k_{\mathrm{inc}}}||$.

In the experiment the relative error is measured against the discrete solution
obtained by solving the scattering problem on the whole domain. In solving the the
scattering problem on the whole domain, the PML is chosen adaptively. These PML
parameters are then fixed and used for all subdomains. Three cases are distinguished.
\begin{enumerate}
\item Schwarz algorithm with two domains (D2): One domain is $\Omega_{1}$ and the
  second domain is the union of $\Omega_{2}$ and $\Omega_{3}$. Thus the layers are
  discretized by finite elements. In Figure~\ref{fig:example1} this corresponds to
  the dark gray lines.
\item Schwarz algorithm with two domains (D2-EUV): One domain is $\Omega_{1}$, the
  second domain is $\Omega_{2}$. $\Omega_{3}$ the layer stack is treated analytically
  and is like a boundary condition for $\Omega_{2}$. That is, if the subproblem on
  $\Omega_{2}$ is solved we iterate internally between $\Omega_{2}$ and $\Omega_{3}$
  and stop if the error is below $10^{-9}$ or after at most $100$ iterations.  In the
  domain decomposition algorithm only the number of iterations between $\Omega_{1}$
  and $\Omega_{2}$ is counted. In Figure~\ref{fig:example1} this corresponds to the
  black lines.
\item Schwarz algorithm with three domains (D3): We are using a multiplicative
  Schwarz algorithm with three subdomains. Within one ``iteration cycle'', we first
  solve for $\Field{E}_{1}$, then for $\Field{E}_{2}$ and finally for
  $\Field{E}_{3}$. In Figure~\ref{fig:example1} this corresponds to the light gray
  lines.
\end{enumerate}
For these three cases the error versus the number of Schwarz iteration cycles is
shown in Figure~\ref{fig:example1} (right).  The experiment is performed for three
different refinement levels, where ``$\times$'' corresponds to the coarsest level
with $5920$ degrees of freedom on the whole domain including the PML, the next finer
level labeled with ``$*$'' is obtained by one uniform refinement of the initial grid
and the finest level labeled with ``$\scriptstyle \Box$'' by two uniform refinements
of the initial grid.
In case (D2-EUV) the error saturates at a level that clearly depends on the refinement of the
interior grid. This behavior can be expected as the number of Fourier coefficients 
that are taken into account to couple the layer-stack analytically in the Schwarz iteration 
is inverse proportional to the mesh-width. In case (D2) and (D3) the error saturates at 
$1e-14$, which is close to machine precision. This surprisingly good convergence behavior
will be further analyzed in a subsequent paper.

\subsection*{Real life EUV mask}

\begin{figure}[tbp]
  \centering
  \psfrag{AS}[cc][cc][0.6][0]{Absorber}
  \psfrag{Mo}[cc][cc][0.6][0]{Mo layer}
  \psfrag{Si}[cc][cc][0.6][0]{Si layer}
  \psfrag{A}[cc][cc][0.6][0]{Air}
  \psfrag{L}[cc][cc][0.6][0]{Line}
  \psfrag{TL}[cc][cc][0.6][0]{Top layer}
  \psfrag{O1}[cc][cc][0.8][0]{$\Omega_{1}$}
  \psfrag{O2}[cc][cc][0.8][0]{$\Omega_{2}$}
  \includegraphics[height=2.5in]{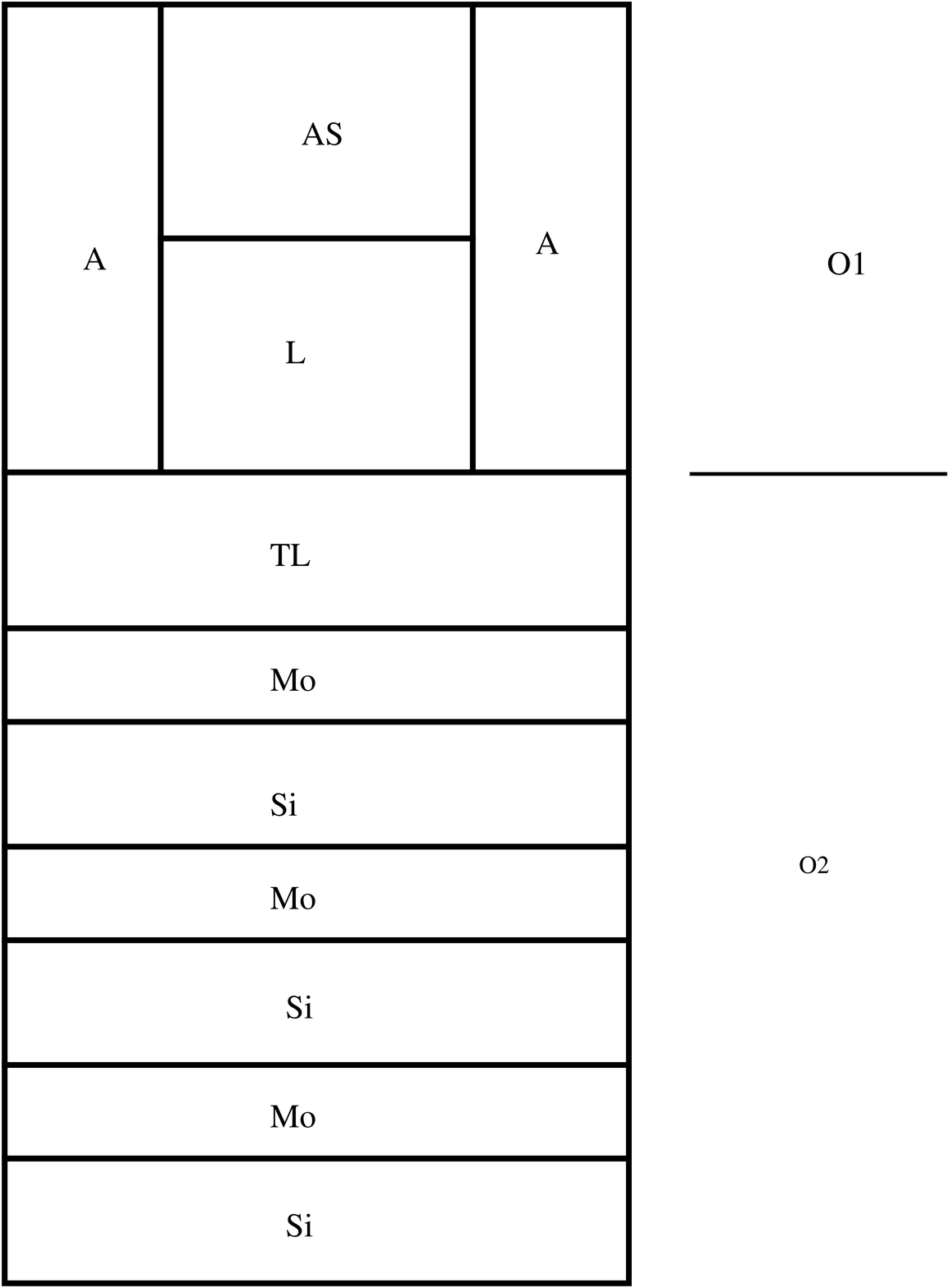}
  \includegraphics[height=2.5in]{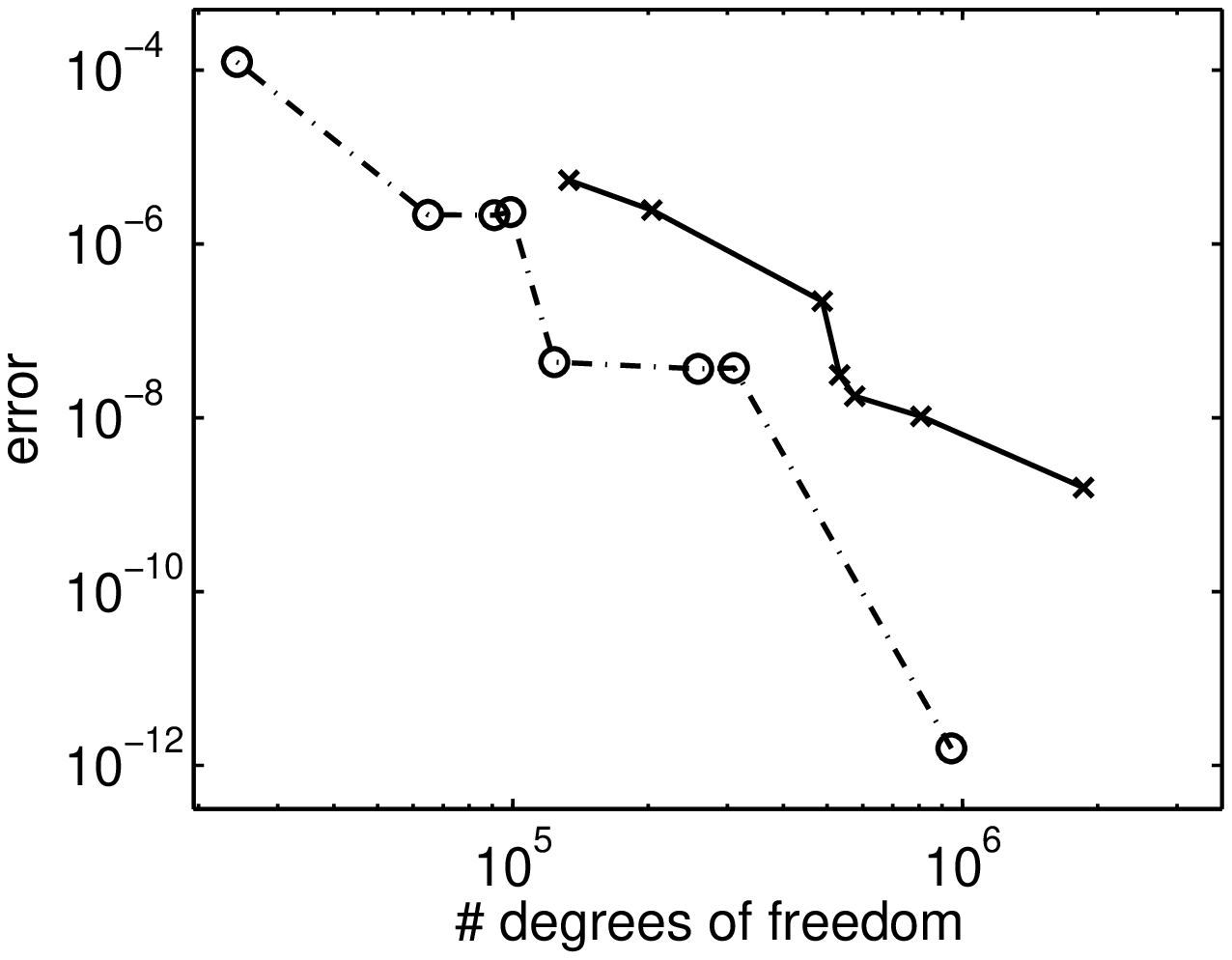} \hspace*{1em} 
  \caption{Left: Sketch of an EUV line mask. 
    Right: Error versus the number of degrees of freedom in finite element mesh.
    The dashed error curve is obtained using the domain decomposition algorithm, decomposing
    the computational domain in two sub-domains (the line and the multi-layer stack) and
    treating the multi-layer stack separately. The solid error curve is obtained 
    discretizing the whole computational domain.  
  }
  \label{fig:erreuvmask}
\end{figure}

A schematic sketch of a more realistic EUV line mask is shown in 
Figure~\ref{fig:erreuvmask}. There only three out of ten MoSi double layers are shown. 
The periodicity $a$ is $40nm$. The line
made of silicon (Si) and the chromium absorber (Cr) have a width of $20nm$ and a
height of $15nm$.  The first silicon layer's height is $10nm$. Each molybdenum
layer (Mo) has a height of $6nm$ and the subsequent silicon layers have a height of
$8nm$. The wavelength is $14nm$. The permeability is $1.0$ everywhere. The
permittivities are $\varepsilon_{Mo}= 1.69+0.016i$, $\varepsilon_{Si}=1.21+0.002i$,
$\varepsilon_{Cr}=1.43+0.24i$ and $\varepsilon_{Air}=1.0$.

Starting from a coarse mesh the grid is pre-refined to have at least $3$, $4$, $5$,
$6$, $7$, $8$, $9$, $16$ and $20$ points per wavelength locally. The solution
obtained with $20$ points per wavelength is taken as a reference solution to measure
the error. 
\\
We use a domain decomposition algorithm and decompose the mask into $\Omega_{1}$
(line, absorber, air) and $\Omega_{2}$ (multilayer-stack). The multilayer-stack is
treated analytically as described in Section~\ref{SubSec:SAEUV}. Additionally we are
using a damping factor of $0.66$ in the domain decomposition algorithm to speed up
convergence. 
\\
The PML is chosen adaptively as described in Section~\ref{SubSec:AdaptedPML}.

Figure~\ref{fig:erreuvmask} shows the error versus the number of degrees of freedom
in the finite element grid including the PML. To obtain the solid line, the multi-layer stack is
discretized using finite elements. Clearly, if the multi-layer stack is not discretized, but
treated analytically and coupled to $\Omega_{1}$ in the domain-decomposition algorithm,
the number of degrees of freedom is reduced drastically. The above calculations where performed
on an {\em{AMD Opeteron PC}} with 16GB of RAM. The arising linear system problems are solved
with the sparse LU method~{\em{PARDISO}},~\cite{Schenk-Gaertner:2004:FFP,Schenk-Gaertner:2004:SUNS}. 

This reduction of the number of degrees of freedom due to the domain decomposition approach, 
allows to compute realistic masks on standard 32--bit computers.

\bibliographystyle{plain} 
\bibliography{group,other}

\end{document}